\documentclass[review,preprint]{elsarticle}

\usepackage[a4paper,total={7in,8in}]{geometry}
\pdfoutput=1
\usepackage{framed,multirow}
\usepackage{graphicx}%
\usepackage{multirow}%
\usepackage{amsmath,amssymb,amsfonts}%
\usepackage{amsthm}%
\usepackage{mathrsfs}%
\usepackage[title]{appendix}%
\usepackage{xcolor}%
\usepackage{textcomp}%
\usepackage{manyfoot}%
\usepackage{booktabs}%
\usepackage{algorithm}%
\usepackage{algorithmicx}%
\usepackage{algpseudocode}%
\usepackage{listings}%
\usepackage{array}
\usepackage{amssymb}
\usepackage{latexsym}
\usepackage{hyperref}
\usepackage{marvosym}

\theoremstyle{thmstyleone}%
\newtheorem{theorem}{Theorem}
\theoremstyle{thmstyletwo}%

\theoremstyle{thmstylethree}%
\usepackage{url}
\usepackage{xcolor}
\definecolor{newcolor}{rgb}{.8,.349,.1}

\begin{document}


\begin{frontmatter}

\title{An efficient unconditional energy stable scheme for the simulation of droplet formation}%

\author[1]{Jinpeng Zhang}

\author[1]{Changjuan Zhang\corref{cor2}}
\cortext[cor2]{Corresponding author: cjzhang@m.scnu.edu.cn}

\author[2,3]{Xiaoping Wang\corref{cor1}}
\cortext[cor1]{Corresponding author: wangxiaoping@cuhk.edu.cn}

\address[1]{South China Research Center for Applied Mathematics and Interdisciplinary Studies, School of Mathematical Sciences, South China Normal University, Guangzhou, 510631, China}
\address[2]{Department of Mathematics, the Hong Kong University of Science and Technology, Hong Kong, China}
\address[3]{School of Science and Engineering, The Chinese University of Hong Kong,
Shenzhen, Guangdong 518172, China $\&$ Shenzhen International Center for Industrial and Applied
Mathematics, Shenzhen Research Institute of Big Data, Guangdong 518172, China}


\begin{abstract}
We have developed an efficient and unconditionally energy-stable method for simulating droplet formation dynamics. Our approach involves a novel time-marching scheme based on the scalar auxiliary variable technique, specifically designed for solving the Cahn-Hilliard-Navier-Stokes phase field model with variable density and viscosity. We have successfully applied this method to simulate droplet formation in scenarios where a Newtonian fluid is injected through a vertical tube into another immiscible Newtonian fluid.
To tackle the challenges posed by nonhomogeneous Dirichlet boundary conditions at the tube entrance, we have introduced additional nonlocal auxiliary variables and associated ordinary differential equations. These additions effectively eliminate the influence of boundary terms. Moreover, we have incorporated stabilization terms into the scheme to enhance its numerical effectiveness. Notably, our resulting scheme is fully decoupled, requiring the solution of only linear systems at each time step. We have also demonstrated the energy decaying property of the scheme, with suitable modifications.
To assess the accuracy and stability of our algorithm, we have conducted extensive numerical simulations. Additionally, we have examined the dynamics of droplet formation and explored the impact of dimensionless parameters on the process.
Overall, our work presents a refined method for simulating droplet formation dynamics, offering improved efficiency, energy stability, and accuracy.
\end{abstract}

\end{frontmatter}


\section{Introduction}
\numberwithin{equation}{section}
Droplet formation is a phenomenon with broad applications in industrial production, such as in pharmaceutics\cite{ref1}, ink-jet printing\cite{ref2}, liquid spraying or atomization\cite{ref3}, DNA micro-arraying\cite{ref4}, etc. The ability to produce stable droplets with controllable sizes is crucial, and research in this area is essential for developing technologies that can achieve this. Numerous experiments\cite{ref5,ref6,ref7,ref8,ref9} have been conducted to advance our understanding of this process.

Numerical simulation is a valuable complement to experimental investigations, and many studies have been conducted to explore this phenomenon. In particular, Wu et al. have reviewed the comprehensive numerical methods used to study fluid dynamics in microfluidic droplet formation\cite{ref10}. These methods can be divided into two classes: sharp interface methods\cite{ref11,ref12,ref13,ref14} and diffuse interface methods\cite{ref15,ref16,ref17}. Xiao, Dianat, and McGuirk presented a sharp interface method that uses a Coupled Level Set/Volume Of Fluid (CLSVOF) technique for interface tracking. Their method accurately predicts droplet formation in low Reynolds number liquid jets as well as the deformation and breakup morphology of a single droplet in uniform air flow at different Weber numbers\cite{ref14}. Matsunaga used a sharp interface method that involves the moving surface mesh particle method to explicitly represent a free-surface boundary\cite{ref18}. This approach enables accurate free-surface tracking and surface tension calculation. They also developed a novel algorithm for determining the particle movement in an ALE fashion by considering a two-dimensional continuity equation. Zhou used finite elements with adaptive meshing in a diffuse-interface framework to simulate the breakup of simple and compound jets in coflowing conditions\cite{ref17}. Liu used a convex splitting scheme for the Cahn-Hilliard equation and a projection type scheme for the momentum equation to accurately simulate the dynamics of droplet formation\cite{ref15}. Although these methods can simulate the dynamics of droplet formation, they require solving complicated equations and can not achieve unconditional energy stability, which makes these methods have high computational costs. Therefore, we aim to develop an easy-to-implement and unconditionally energy-stable numerical method to efficiently solve the droplet formation problem.

The diffuse interface method is a powerful numerical technique for handling topological changes in interfaces, which is important in applications such as droplet formation. Our model is based on a phase field approach that consists of a coupled system of the Cahn-Hilliard and Navier-Stokes equations. For the Cahn-Hilliard-Navier-Stokes (CH-NS) model, several highly successful schemes have been developed, including those cited in references \cite{ref19,ref21,ref22}. These schemes exhibit fully-decoupled characteristics and unconditional energy stability, making them ideal for simulating interfacial phenomena. In particular, \cite{ref22} introduced a new fully decoupled numerical scheme known as the Decoupled Scalar Auxiliary Variable (DSAV) method. The DSAV method achieves unconditional energy stability by explicitly discretizing nonlinear coupling terms. It combines the penalty method of the Navier-Stokes equations with the Strang operator splitting method and introduces nonlocal auxiliary variables and their associated ordinary differential equations. These variables are used to handle various coupled nonlinear terms, including advection and surface tensions.

In this paper, we generalize the DSAV method to the droplet formation problem which is modeled by Cahn-Hilliard-Navier-Stokes equations with a non-homogeneous Dirichlet boundary. To design an unconditionally energy stable numerical scheme for the problem, we introduce more nonlocal auxiliary variables and associated ordinary differential equations to eliminate the boundary terms. Moreover, we add some stabilization terms to enhance the effectiveness of our algorithm and propose a new modified penalty method that also decouples the computation of pressure from the momentum equation, enabling us to demonstrate the unconditional energy stability of our numerical scheme. After overcoming these challenges, our scheme still retains the advantages of the DSAV method, such as (i) it can solve the Cahn-Hilliard equation with constant coefficients; (ii) it is explicit for all nonlinear coupling terms in the fluid momentum equation, aiming to minimize computational costs as much as possible; and (iii) it can maintain linearity, decoupling format, and ensures unconditional energy stability.

We then investigate the accuracy, energy stability, and effectiveness of the proposed scheme numerically. We perform several numerical simulations to confirm the convergence rate and energy stability of the method. We also examine how the dynamics of droplet formation depend on various physical parameters of the system and compare our numerical results with physical experiments. The simulation results demonstrate that the process of drop formation can be reasonably predicted by the phase field model we used.

The remainder of the paper is organized as follows. In Section 2, we provide a concise description of the mathematical formulation of the problem and derive the energy law for the PDE system. We then establish the modified model and prove the energy law for the modified version. In Section 3, we construct the numerical scheme and offer rigorous proof of its unconditional energy stability. Additionally, we provide information on its solvability and detailed implementation process. In Section 4, we present several numerical examples that demonstrate the stability, accuracy, and efficacy of our proposed scheme. We also examine the effects of various dimensionless parameters on the maximum radius of the droplet.

\section{The phase field model}

\label{sec1}
\subsection{Governing equations and boundary conditions}

\begin{figure}[!t] 
\centering 
\includegraphics[width=80mm]{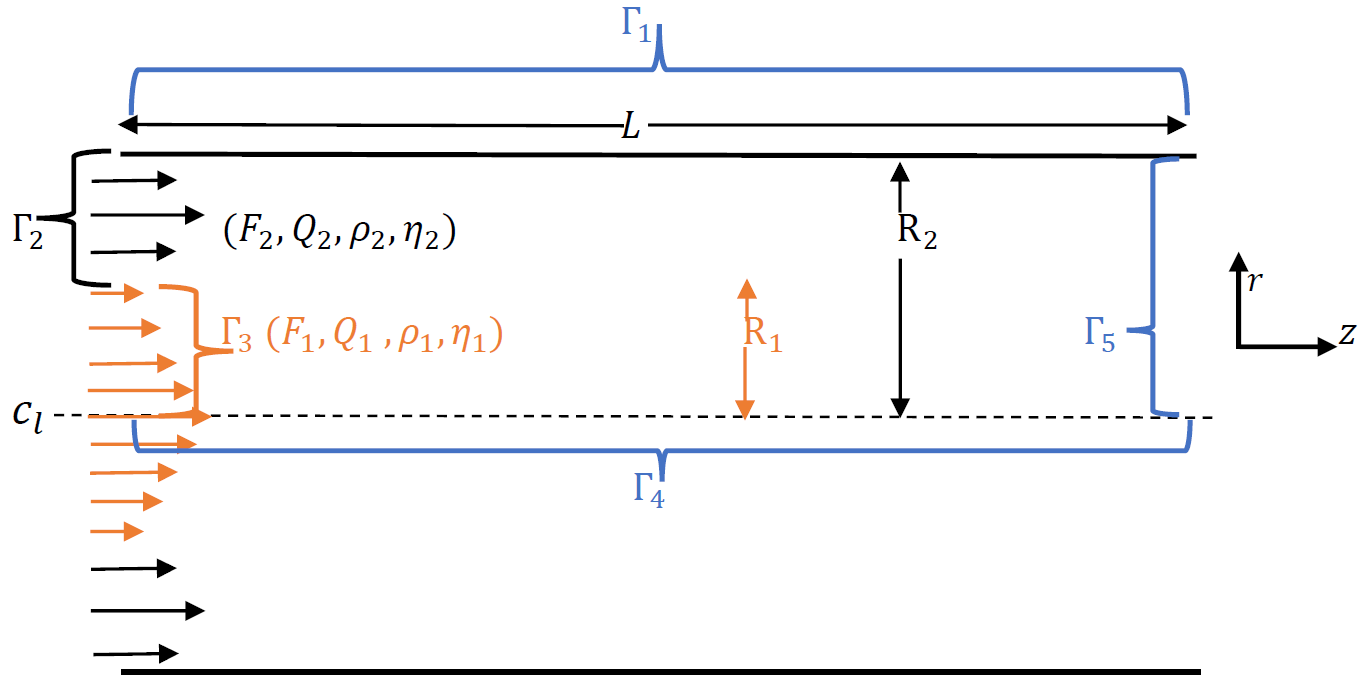} 
\caption{Schematic diagram of droplet formation in another coflowing fluid.} 
\label{fig:Fig 1}
\end{figure}

Consider the injection an incompressible Newtonian fluid $F_{1}$ with density $\rho_{1}$ and viscosity $\eta_{1}$ into another immiscible, incompressible, coflowing Newtonian fluid $F_{2}$ with density $\rho_{2}$ and viscosity $\eta_{2}$ through a vertical capillary tube with a radius of $R_{1}$. The outer fluid is contained in a coaxial cylindrical tube with a radius of $R_{2}$. These two fluids are injected into the tube at constant flow rates of $Q_{1}$ and $Q_{2}$ respectively. Fig.\ref{fig:Fig 1} shows a schematic diagram. We use a cylindrical coordinate system $\{r,z,\theta\}$ with its origin located at the intersection of the centerline $c_{l}$ and the inflow boundary $z = 0$. Here $\{r,z,\theta\}$ represent the radial coordinate, axial coordinate, and azimuthal angle respectively. We assume radial symmetry, implying that all variables are independent of the azimuthal angle $\theta$.  The computational domain $\Omega$ is limited to the upper half of the tube. The domain is bounded by solid wall denoted as $\Gamma_{1}$, the inlet of the tube denoted as $\partial\Omega_{L}=\Gamma_{2}\cup\Gamma_{3}$ , the outlet of the tube specified by $\partial\Omega_{R}=\Gamma_{5}$, and the center line identified as $\Gamma_{4}$, as illustrated in Fig.\ref{fig:Fig 1} 

The governing equations of the CH-NS model with variable density and viscosity in the cylindrical coordinate system read as
\begin{align}
&\phi_{t}+\mathbf{u}\cdot\nabla\phi=\frac{\mathcal{M}}{r}\nabla\cdot r\nabla\mu \label{21}, \\
&\mu=-\frac{\mathcal{J}}{r}\nabla\cdot\left(r\nabla\phi\right)-e\phi+k\phi^{3} \label{22},\\
&\rho(\phi)\left(\mathbf{u}_{t}+(\mathbf{u} \cdot \nabla) \mathbf{u}\right)+\mathbf{J}(\mu)\cdot\nabla\mathbf{u}-\frac{1}{r}\nabla \cdot(\eta(\phi)rD(\mathbf{u}))+\nabla p-\mu \nabla \phi +\left(0,2\eta(\phi)\frac{v_{r}}{r^{2}}\right)=0 \label{23}, \\
&\frac{1}{r}\nabla \cdot (r\mathbf{u})=0 \label{24}.
\end{align}
Here $\nabla = \left(\frac{\partial }{\partial z},\frac{\partial }{\partial r} \right)$, $\phi$ is a phase-field variable, $\mathbf{u}=(v_{z},v_{r})$ is the fluid velocity where $v_{z}$ and $v_{r}$ represent the axial and radial components respectively, $\mathcal{M}$ is the phenomenological mobility coefficient, $\mu$ is the chemical potential, $D(\mathbf{u})=\nabla\mathbf{u}+(\nabla\mathbf{u})^{T}$ denotes the viscous part of the stress tensor, $\mathbf{J}(\mu)=\mathcal{M}\nabla\mu\frac{\rho_{1}-\rho_{2}}{2}$, $p$ is the pressure. Parameters $\mathcal{J}$, $e$, $k$ in \eqref{22} are related with the interface thickness $\xi=\sqrt{\frac{\mathcal{J}}{e}}$, the interfacial tension $\gamma=\frac{2\sqrt{2}e^{2}\xi}{3k}$ and $\phi_{\pm}=\pm\sqrt{\frac{e}{k}}$(=$\pm 1$ in our case). $(\mathbf{u} \cdot \nabla) \mathbf{u}$ and $\nabla\cdot(\mathbf{u}\phi)$ are two advection terms. Density and viscosity are assumed as an interpolation function of $\phi$,
\begin{equation}
\rho(\phi)=\rho_{1}\frac{1-\phi}{2}+\rho_{2}\frac{1+\phi}{2}, \quad
\eta(\phi)=\eta_{1}\frac{1-\phi}{2}+\eta_{2}\frac{1+\phi}{2}. \label{25}
\end{equation}

No-slip and no-penetration conditions are imposed along the solid walls of the tubes, 
\begin{equation}
\mathbf{u}(z, r, t)=0,\quad\frac{\partial\mu}{\partial r}(z, r, t)=0,\quad\frac{\partial\phi}{\partial r}(z, r, t)=0, \quad \text{on} \quad\Gamma_{1}.
\end{equation}
At the inlet of the tube, the inflow conditions for $\mathbf{u}$ at $z = 0$ are 
\begin{equation}
v_r(z, r, t)=0, \quad v_z(z, r, t)=2 \frac{Q_2}{\pi R_2^2}\left[\frac{1-\left(\frac{r}{R_2}\right)^2+\frac{1-\left(R_1 / R_2\right)^2}{\ln \left(R_2 / R_1\right)} \ln \left(\frac{r}{R_2}\right)}{1-\left(\frac{R_1}{R_2}\right)^4-\frac{\left(1-\left(R_1 / R_2\right)^2\right)^2}{\ln \left(R_2 / R_1\right)}}\right],\quad \text{on} \quad \Gamma_{2},\label{V_b1}
\end{equation}
\begin{equation}
v_r(z, r, t)=0, \quad v_z(z, r, t)=2 \frac{Q_1}{\pi R_1^2}\left[1-\left(\frac{r}{R_1}\right)^2\right],\quad \text{on}\quad \Gamma_{3}, \label{V_b2}
\end{equation}
The velocity profiles \eqref{V_b1}-\eqref{V_b2} are also used by Suryo and Basaran\cite{Suryo}. The inflow boundary conditions for $\phi$ and $\mu$ are
\begin{equation}
 \phi(z, r, t)=1,\quad \mu(z, r, t)=0,\quad \text{on}\quad\Gamma_{2},
\end{equation}
\begin{equation}
\phi(z, r, t)=-1,\quad \mu(z, r, t)=0,\quad \text{on} \quad\Gamma_{3}.
\end{equation}
Along the central line $c_{l}$, we use symmetric boundary conditions,
\begin{equation}
v_r(z, r, t)=0, \quad \frac{\partial v_z}{\partial r}(z, r, t)=0,\quad\frac{\partial\phi}{\partial r}(z, r, t)=0,\quad \frac{\partial\mu}{\partial r}(z, r, t)=0,\quad \text{on} \quad\Gamma_{4}.
\end{equation}
Assuming that the length of the outer tube is long enough, so that the outflow condition  will not affect the droplet formation process, and at the outlet boundary $z = L$, we can assume
\begin{equation}
v_{r}(z, r, t)=0,\quad\frac{\partial v_{z}}{\partial z}(z, r, t)=0,\quad\frac{\partial\mu}{\partial z}(z, r, t)=0,\quad\frac{\partial\phi}{\partial z}(z, r, t)=0,\quad \text{on} \quad\Gamma_{5}.
\end{equation}
The initial conditions read as,
\begin{equation}
(\phi,\mathbf{u},\mu)|_{t=0}=(1,\mathbf{0},0). \label{26}
\end{equation}
We now introduce the following characteristic scales,
$$
\begin{aligned}
l_{c}=R_{1},\quad v_{c}=\frac{Q_{1}}{\pi R_{1}^{2}},\quad\rho_{c}=\rho_{1},\quad\eta_{c}=\eta_{1},
\end{aligned}
$$
The dimensionless counterparts of the system \eqref{21}-\eqref{24} are as follows, where we have utilized the same notations for the dimensionless variables.
\begin{align}
&\phi_{t}+\mathbf{u}\cdot\nabla\phi=\frac{\mathcal{L}_{d}}{r}\nabla\cdot (r\nabla\mu), \label{27} \\
&\mu=-\frac{\epsilon}{r}\nabla\cdot\left(r\nabla\phi\right)+\frac{1}{\epsilon} f(\phi), \label{28}\\
&\mathcal{R}e\rho(\phi)\left(\mathbf{u}_{t}+(\mathbf{u} \cdot \nabla) \mathbf{u}\right)+\mathbf{J}(\mu)\cdot\nabla\mathbf{u}-\frac{1}{r}\nabla \cdot(\eta(\phi)rD(\mathbf{u}))+\nabla p-\mathcal{B}\mu \nabla \phi+\left(0,2\eta(\phi)\frac{v_{r}}{r^{2}}\right)=0,\label{29} \\
&\frac{1}{r}\nabla \cdot (r\mathbf{u})=0, \label{210}\\
&\rho(\phi)=\frac{1-\phi}{2}+\lambda_{\rho}\frac{1+\phi}{2}, 
\eta(\phi)=\frac{1-\phi}{2}+\lambda_{\eta}\frac{1+\phi}{2}.\label{211}
\end{align}
The boundary conditions of this system read as
\begin{equation}
    \begin{split}
&\mathbf{u}(z, r, t)=0,\quad\frac{\partial\mu}{\partial r}(z, r, t)=0,\quad\frac{\partial\phi}{\partial r}(z, r, t)=0,\quad \text{on} \quad\Gamma_{1}, \\
&v_r(z, r, t)=0, \quad v_z(z, r, t)=2 \frac{Q_r}{a^2}\left[\frac{1-\left(\frac{r}{a}\right)^2+\frac{1-\left(1 /a\right)^2}{\ln \left(a\right)} \ln \left(\frac{r}{a}\right)}{1-\left(\frac{1}{a}\right)^4-\frac{\left(1-\left(1 / a\right)^2\right)^2}{\ln \left(a\right)}}\right], \quad\phi(z, r, t)=1,\quad
\mu(z, r, t)=0, \quad \text{on} \quad\Gamma_{2}, \\
&v_r(z, r, t)=0, \quad v_z(z, r, t)=2\left[1-r^2\right],\quad
\phi(z, r, t)=-1,\quad \mu(z, r, t)=0,\quad \text{on} \quad\Gamma_{3},\\
&v_r(z, r, t)=0, \quad \frac{\partial v_z}{\partial r}(z, r, t)=0,\quad\frac{\partial\phi}{\partial r}(z, r, t)=0,\quad
\frac{\partial\mu}{\partial r}(z, r, t)=0,\quad \text{on} \quad\Gamma_{4},\\
&v_{r}(z, r, t)=0,\quad\frac{\partial v_{z}}{\partial z}(z, r, t)=0,\quad\frac{\partial\phi}{\partial z}(z, r, t)=0,\quad
\frac{\partial\mu}{\partial z}(z, r, t)=0,\quad \text{on} \quad\Gamma_{5}.
    \end{split}\label{212}
\end{equation}
Here $f(\phi)=\phi^{3}-\phi$. $\mathbf{J}(\mu)=\mathcal{L}_{d}\mathcal{R}e\nabla\mu(1-\lambda_{\rho})/2$ is the dimensionless form.
The dimensionless parameter Reynolds number measures the relative importance of inertial force to viscous force, given by  $\mathcal{R}e=\rho_{c}v_{c}l_{c}/\eta_{c}$. Capillary number measures the relative importance of the viscous force to surface tension force, given by $\mathcal{C}a=\eta_{c}v_{c}/\gamma$. Cahn number is the ratio between interface thickness and length scale, given by $\epsilon=\xi/l_{c}$. $\mathcal{L}_{d}=3\mathcal{M}\gamma/(2\sqrt{2}v_{c}\xi^{2})$ is the diffusion coefficient. $\lambda_{\rho}=\rho_{2}/\rho_{1}$ and $\lambda_{\eta}=\eta_{2}/\eta_{1}$ are the density ratio and viscosity ratio respectively. $Q_{r}=Q_{2}/Q_{1}$ is the ratio of the flow rate of the outer fluid to that of the inner fluid.  $a=R_{2}/R_{1}$ is the ratio of the radius.\\

The energy we proposed is different from the energy in \cite{ref22} because of the inflow condition. We derive the law of energy dissipation for the system \eqref{27}-\eqref{210} as follows\\
\begin{theorem}
The solutions of the system \eqref{27}-\eqref{210},with boundary conditions \eqref{212}, \textcolor{red}{satisfy}
\begin{align}
\frac{d}{dt}\mathbf{E}_{O}=-\frac{1}{2}\int_{\Omega}\eta r|D(\mathbf{u})|^{2}drdz-\int_{\Omega}2\eta\frac{v_{r}^{2}}{r}drdz-\mathcal{B}\mathcal{L}_{d}\int_{\Omega}r|\nabla\mu|^{2}drdz \leq 0,
\end{align}
where, 
\begin{equation}
    \begin{split}
&\mathbf{E}_{O}=\frac{\mathcal{R}e}{2}\int_{\Omega}r|\sigma\mathbf{u}|^{2}drdz+\frac{\mathcal{B}\epsilon}{2}\int_{\Omega}r|\nabla\phi|^{2}drdz+\frac{\mathcal{B}}{\epsilon}\int_{\Omega}rF(\phi)drdz \\ 
&+\int_{0}^{t}\Big\{-\int_{\partial\Omega_{L}}rpv_{z}ds-\frac{\mathcal{R}e}{2}\int_{\partial\Omega_{L}}r\rho v_{z}^{3}ds+\frac{\mathcal{R}e}{2}\int_{\partial\Omega_{R}}r\rho v_{z}^{3}ds+2\int_{\partial\Omega_{L}}r\eta v_{z}\frac{\partial v_{z}}{\partial z}ds\\ 
&+\int_{\partial\Omega_{R}}rpv_{z}ds-\mathcal{L}_{d}\mathcal{R}e\frac{1-\lambda_{\rho}}{4}\int_{\partial\Omega_{L}}\frac{\partial\mu}{\partial z}v_{z}^{2}ds\Big\} d\tau. 
    \end{split}\label{217}
\end{equation}
\end{theorem}
\noindent\textbf{Proof:} By multiplying $r\mathcal{R}e\frac{\lambda_{\rho}-1}{4}$ to \eqref{27}, we derive
\begin{equation}
r\frac{\mathcal{R}e}{2}\rho_{t}+\frac{\mathcal{R}e}{2}\nabla\cdot\left(r\rho\mathbf{u}\right)+\frac{1}{2}\nabla\cdot r\mathbf{J}(\mu)=0.
\end{equation}
Futhermore, by multiplying $\mathbf{u}$ to this equation,
\begin{align}
r\frac{\mathcal{R}e}{2}\rho_{t}\mathbf{u}+\frac{\mathcal{R}e}{2}\nabla\cdot\left(r\rho\mathbf{u}\right)\mathbf{u}+\frac{1}{2}\nabla\cdot r\mathbf{J}(\mu)\mathbf{u}=0. \label{218}
\end{align}
By combining \eqref{29} and \eqref{218},
\begin{equation}
    \begin{split}
&\mathcal{R}e\int_{\Omega}r\sigma\mathbf{u}\cdot\left(\sigma\mathbf{u}\right)_{t}drdz+\mathcal{R}e\int_{\Omega}\left(\frac{1}{2}\nabla\cdot\left(r\rho\mathbf{u}\right)\mathbf{u}\cdot\mathbf{u}+r\rho(\mathbf{u} \cdot \nabla) \mathbf{u}\cdot\mathbf{u}\right)drdz\\
&-\int_{\Omega}\nabla\cdot(r\eta D(\mathbf{u}))\cdot\mathbf{u}drdz+\int_{\Omega}\frac{1}{2}\nabla\cdot r\mathbf{J}(\mu)\mathbf{u}\cdot\mathbf{u}+r\mathbf{J}(\mu)\cdot\nabla\mathbf{u}\cdot\mathbf{u}drdz\\
&\int_{\Omega}r\nabla p\cdot\mathbf{u}drdz-\mathcal{B}\int_{\Omega}r\mu \nabla \phi\cdot\mathbf{u}drdz+\int_{\Omega}2\eta\frac{v_{r}^{2}}{r}drdz=0,
    \end{split}
\end{equation}
where $\sigma=\sqrt{\rho}$. By taking the $L^{2}$ inner product of this equation with $\mathbf{u}$, using integration by parts, and applying the divergence-free condition \eqref{210}, we obtain
\begin{equation}
    \begin{split}
&-\frac{\mathcal{R}e}{2}\int_{\partial\Omega_{L}}r\rho v_{z}^{3}ds+\frac{\mathcal{R}e}{2}\int_{\partial\Omega_{R}}r\rho v_{z}^{3}ds+\frac{\mathcal{R}e}{2}\frac{\partial}{\partial t}\int_{\Omega}r|\sigma\mathbf{u}|^{2}drdz+\frac{1}{2}\int_{\Omega}r\eta |D(\mathbf{u})|^{2}drdz \\ 
&-\int_{\partial\Omega_{L}}rpv_{z}ds+\int_{\partial\Omega_{R}}rv_{z}p ds+2\int_{\partial \Omega_{L}}r\eta v_{z}\frac{\partial v_{z}}{\partial z}drdz-\mathcal{B}\int_{\Omega}r\mu \nabla \phi\cdot\mathbf{u}drdz  \\
&-\mathcal{L}_{d}\mathcal{R}e\frac{1-\lambda_{\rho}}{2}\int_{\partial\Omega_{L}}\frac{\partial\mu}{\partial z}v_{z}^{2}ds+\int_{\Omega}2\eta\frac{v_{r}^{2}}{r}drdz=0.
    \end{split}\label{219}
\end{equation}
By taking an inner product of \eqref{27} with $r\mathcal{B}\mu$, and performing integration by parts, we derive

\begin{align}
\mathcal{B}\int_{\Omega}r\phi_{t}\mu drdz+\mathcal{B}\int_{\Omega}r(\mathbf{u}\cdot\nabla\phi)\mu drdz=-\mathcal{B}\mathcal{L}_{d}\int_{\Omega}r|\nabla\mu|^{2}drdz. \label{220}
\end{align}
Then, we take the inner product of \eqref{28} with $r\mathcal{B}\phi_{t}$ in the $L^{2}$ space, 
\begin{equation}
\mathcal{B}\int_{\Omega}r\phi_{t}\mu drdz=-\mathcal{B}\epsilon\int_{\Omega}\nabla\cdot\left(r\nabla\phi\right)\phi_{t}drdz+\frac{\mathcal{B}}{\epsilon}\int_{\Omega}rf(\phi)\phi_{t}drdz,
\end{equation}
and use integration by parts to get
\begin{align}
\mathcal{B}\int_{\Omega}r\phi_{t}\mu drdz=\frac{\mathcal{B}\epsilon}{2}\frac{\partial }{\partial t}\int_{\Omega}r|\nabla\phi|^{2}drdz+\frac{\mathcal{B}}{\epsilon}\frac{\partial }{\partial t}\int_{\Omega}rF(\phi)drdz.  \label{221}
\end{align}
From  \eqref{219}, \eqref{220} and \eqref{221}, we derive
\begin{equation}
    \begin{split}
    &\frac{\mathcal{R}e}{2}\frac{\partial}{\partial t}\int_{\Omega}r|\sigma\mathbf{u}|^{2}drdz+\frac{\mathcal{B}\epsilon}{2}\frac{\partial }{\partial t}\int_{\Omega}r|\nabla\phi|^{2}drdz+\frac{\mathcal{B}}{\epsilon}\frac{\partial }{\partial t}\int_{\Omega}rF(\phi)drdz+2\int_{\partial\Omega_{L}}r\eta v_{z}\frac{\partial v_{z}}{\partial z}ds\\
    &-\int_{\partial\Omega_{L}}rpv_{z}ds+\int_{\partial\Omega_{R}}rpv_{z}ds-\frac{\mathcal{R}e}{2}\int_{\partial\Omega_{L}}r\rho v_{z}^{3}ds-\mathcal{L}_{d}\mathcal{R}e\frac{1-\lambda_{\rho}}{4}\int_{\partial\Omega_{L}}\frac{\partial\mu}{\partial z}v_{z}^{2}ds+\frac{\mathcal{R}e}{2}\int_{\partial\Omega_{R}}r\rho v_{z}^{3}ds\\
    =&-\frac{1}{2}\int_{\Omega}\eta r|D(\mathbf{u})|^{2}drdz-\int_{\Omega}2\eta\frac{v_{r}^{2}}{r}drdz-\mathcal{B}\mathcal{L}_{d}\int_{\Omega}r|\nabla\mu|^{2}drdz \leq 0.
    \end{split}\label{222}
\end{equation}
We propose a definition
\begin{equation}
    \begin{split}
    \mathbf{E}_{O}=&\frac{\mathcal{R}e}{2}\int_{\Omega}r|\sigma\mathbf{u}|^{2}drdz+\frac{\mathcal{B}\epsilon}{2}\int_{\Omega}r|\nabla\phi|^{2}drdz+\frac{\mathcal{B}}{\epsilon}\int_{\Omega}rF(\phi)drdz+\int_{0}^{t}\Big\{-\int_{\partial\Omega_{L}}rpv_{z}ds+\int_{\partial\Omega_{R}}rpv_{z}ds\\
    &-\frac{\mathcal{R}e}{2}\int_{\partial\Omega_{L}}r\rho v_{z}^{3}ds+\frac{\mathcal{R}e}{2}\int_{\partial\Omega_{R}}r\rho v_{z}^{3}ds+2\int_{\partial\Omega_{L}}r\eta v_{z}\frac{\partial v_{z}}{\partial z}ds-\mathcal{L}_{d}\mathcal{R}e\frac{1-\lambda_{\rho}}{4}\int_{\partial\Omega_{L}}\frac{\partial\mu}{\partial z}v_{z}^{2}ds \Big\} d\tau.
    \end{split}\label{223}
\end{equation}
Then, the law of energy dissipation of the CH-NS system \eqref{27}-\eqref{210} is given as follows,
\begin{equation}
    \begin{split}
    \frac{d}{dt}\mathbf{E}_{O}=-\frac{1}{2}\int_{\Omega}\eta r|D(\mathbf{u})|^{2}drdz-\int_{\Omega}2\eta\frac{v_{r}^{2}}{r}drdz-\mathcal{B}\mathcal{L}_{d}\int_{\Omega}r|\nabla\mu|^{2}drdz \leq 0.
    \end{split}
\end{equation}

\subsection{Modified model and its energy law}
Due to the influence of nonlinear coupling terms and non-homogeneous Dirichlet boundary conditions, if we directly discretize the original system \eqref{27}-\eqref{210}, it is difficult to develop a linear, fully-decoupled and unconditionally energy-stable numerical scheme. In the derivation of energy law, we can see many non-zero boundary terms, which is the main difficulty in developing an effective numerical method. \cite{ref22} proposed a novel fully-decoupled numerical technique that can achieve unconditional energy stability while explicitly discretizing nonlinear coupling items, such as $f(\phi)$ in Cahn-Hilliard equation, $(\mathbf{u}\cdot\nabla)\mathbf{u}$ in the momentum equation, and $\rho\mathbf{u}_{t}$,  $\rho_{t}\mathbf{u}$ associated with time derivatives. However, it can not be directly applied to the problem of droplet formation because of the influence of the non-zero boundary terms. We will extend this method to solve the droplet formation problem, the main difficulty is how to handle the non-homogeneous boundary conditions.

To address this challenge, we introduce nonlocal variables and design special ODEs for them. Firstly, to eliminate the effect of non-zero boundary terms, we introduce a nonlocal variable $\boldsymbol{K}$ that reads as
\begin{align}
\boldsymbol{K}=\sqrt{\int_{0}^{t}-\mathcal{K}d\tau+\boldsymbol{G}},\label{bK}
\end{align}
where $\boldsymbol{G}$ is a positive constant and 

\begin{align}
\mathcal{K}=&\frac{\mathcal{R}e}{2}\int_{\partial\Omega_{L}}r\rho v_{z}^{3}d \boldsymbol{s}-\frac{\mathcal{R}e}{2}\int_{\partial\Omega_{R}}r\rho v_{z}^{3}d \boldsymbol{s}+\int_{\partial\Omega_{L}}rpv_{z}ds-\int_{\partial\Omega_{R}}rpv_{z}ds-2\int_{\partial\Omega_{L}}r\eta v_{z}\frac{\partial v_{z}}{\partial z}ds\\
&+\mathcal{L}_{d}\mathcal{R}e\frac{1-\lambda_{\rho}}{4}\int_{\partial\Omega_{L}}\frac{\partial\mu}{\partial z}v_{z}^{2}ds.\nonumber \label{K}
\end{align}

The variable $\boldsymbol{K}$ is used to "quadratize" the non-zero boundary terms. Since our focus is exclusively on the droplet formation process occurring within a finite time, $\int_{0}^{t}-\mathcal{K}d\tau$ can be bounded below. The constant $\boldsymbol{G}$ is used to ensure that the radicand is positive. This quadratization approach to handling the non-zero boundary terms using a nonlocal auxiliary variable is similar to the so-called SAV method \cite{ref24, ref25} which is an efficient method to linearize the nonlinear terms induced by the energy potentials.

Secondly, to process coupled nonlinear terms, such as advection and surface tensions, and eliminate the effect of boundary terms produced by these coupled nonlinear terms, we introduce three nonlocal variables $\boldsymbol{Q}$, $\boldsymbol{R}$ and $\boldsymbol{T}$ and design special ODEs for each of them.
\begin{align}
&\left\{\begin{array}{l}
\boldsymbol{Q}_{t}=\alpha\int_{\Omega}\left(r\mathbf{u}\cdot\nabla\phi \mu-r\mathbf{u}\cdot\nabla\phi \mu\right) drdz,\\
\boldsymbol{Q}|_{t=0}=1.
\end{array}\right. \\
&\left\{\begin{array}{l}
\boldsymbol{R}_{t}=\alpha\bigg(\int_{\Omega}\left(\mathcal{R}e r\rho(\phi)(\mathbf{u} \cdot \nabla) \mathbf{u} \cdot \mathbf{u}+\frac{\mathcal{R}e}{2}\nabla \cdot(r\rho(\phi)\mathbf{u})\mathbf{u} \cdot \mathbf{u}+r\nabla p\cdot \mathbf{u}\right)drdz\\
-\int_{\Omega}\big(\nabla \cdot(r\eta(\phi)D(\mathbf{u}))\big)\cdot\mathbf{u}drdz-\frac{1}{2}\int_{\Omega}r\eta|D(\mathbf{u})|^{2}drdz+\\
\int_{\Omega}\frac{1}{2}\nabla\cdot r\mathbf{J}(\mu)\mathbf{u}\cdot\mathbf{u}+r\mathbf{J}(\mu)\cdot\nabla\mathbf{u}\cdot\mathbf{u}drdz+\frac{\boldsymbol{K}\mathcal{K}}{\sqrt{\int_{0}^{t}-\mathcal{K}d\tau+\boldsymbol{G}}}\bigg),\\
\boldsymbol{R}|_{t=0}=1.
\end{array}\right. \\
&\left\{\begin{array}{l}
\boldsymbol{T}_{t}=\alpha\int_{\Omega}\left(\frac{\partial (rv_{r})}{\partial r}+\frac{\partial (rv_{z})}{\partial z}\right)p drdz,\\
\boldsymbol{T}|_{t=0}=1.
\end{array}\right. 
\end{align}\\
It is easy to see that these three ODEs are equivalent to $\boldsymbol{Q}_{t}=0$, $\boldsymbol{R}_{t}=0$ and $\boldsymbol{T}_{t}=0$. Thus $\boldsymbol{Q}=1$, $\boldsymbol{R}=1$ and $\boldsymbol{T}=1$ are their exact solutions. $\alpha$ is the stabilization parameter, usually taken as a very small parameter to ensure that the number at the right end of the equation is approximately zero after numerical discretization. Several other studies share a common approach of incorporating stabilization parameters to improve algorithmic properties. For instance, Ju et al. \cite{ju2022stabilized,ju2022generalized} presented the stabilized exponential-SAV method to preserve the maximum bound principle by introducing a stabilization parameter.

Thirdly, we add a stabilization term $\frac{(1-\mathbf{R})}{r}\nabla \cdot(r\eta(\phi)D(\mathbf{u}))$ in the momentum equation to enhance the stability of our scheme. In the next section, we will demonstrate how to improve numerical stability through this term. Due to $\mathbf{R}=1$, it is easy to know the exact value of this term is zero, so it will not change the equation \eqref{29}.\par
Fourthly, we introduce another nonlocal scalar variable $\boldsymbol{U}$, which is defined as
\[
\boldsymbol{U}=\sqrt{\int_{\Omega}\left(rF(\phi)-\frac{rs}{2}\phi^{2}\right) drdz+B}=\sqrt{\int_{\Omega}\left(r\frac{(\phi^{2}-1)^{2}}{4}-\frac{rs}{2}\phi^{2}\right) drdz+B},
\]
where $s$ and $B$ are two positive constants. This variable $\boldsymbol{U}$ is used to "quadratize" the nonlinear double-well potential. We extract the quadratic term $\frac{s}{2}\phi^{2}(s\sim O(1))$ from the double-well potential $F(\phi)$ which can help to maintain the $H^{1}$ stability of $\phi$ (see \cite{ref28}). The constant $B$ is used to ensure the radicand positive since $F(\phi)-\frac{s}{2}\phi^{2}$ is always bounded from below.\par
By using these new variables, we can rewrite the CH-NS system \eqref{27}-\eqref{210} to the equivalent form which is called the modified model:
\begin{align}
&\phi_{t}+\boldsymbol{Q}\mathbf{u}\cdot\nabla\phi=\frac{1}{r}\mathcal{L}_{d}\nabla\cdot r\nabla\mu, \label{224}\\
&\mu=-\frac{\epsilon}{r} \nabla\cdot\left(r\nabla\phi\right)+\frac{s}{\epsilon}\phi+\frac{1}{\epsilon}\boldsymbol{H}(\phi)\boldsymbol{U}, \label{225}\\
&\mathcal{R}e\left(\frac{1}{2}\rho_{t}(\phi)\mathbf{u}+\rho(\phi)\mathbf{u}_{t}+\boldsymbol{R}     \rho(\phi)(\mathbf{u} \cdot \nabla) \mathbf{u}+\frac{\boldsymbol{R}}{2r}\nabla \cdot(r\rho(\phi)\mathbf{u})\mathbf{u}\right)-\frac{1}{r}\nabla \cdot(r\eta(\phi)D(\mathbf{u}))+\mathbf{R}\nabla p \nonumber\\
&-\boldsymbol{Q}\mathcal{B}\mu \nabla \phi+\left(0,2\eta(\phi)\frac{v_{r}}{r^{2}} \right)+\frac{(1-\mathbf{R})}{r}\nabla \cdot(r\eta(\phi)D(\mathbf{u}))+\frac{\boldsymbol{R}}{2}\nabla\cdot r\mathbf{J}(\mu)\mathbf{u}+r\boldsymbol{R}\mathbf{J}(\mu)\cdot\nabla\mathbf{u}=0,\label{226}
\end{align}
\begin{align}
&\frac{1}{r}\nabla\cdot\left(r\mathbf{u} \right)=0,\label{232}\\ 
&\boldsymbol{U}_{t}=\frac{1}{2}\int_{\Omega}r\boldsymbol{H}(\phi)\phi_{t}drdz, \label{227}\\
&\boldsymbol{Q}_{t}=\alpha\int_{\Omega}\left(r\mathbf{u}\cdot\nabla\phi \mu-r\mathbf{u}\cdot\nabla\phi \mu\right) drdz, \label{228}\\
&\boldsymbol{R}_{t}=\alpha\bigg(\int_{\Omega}\left(\mathcal{R}e r\rho(\phi)(\mathbf{u} \cdot \nabla) \mathbf{u} \cdot \mathbf{u}+\frac{\mathcal{R}e}{2}\nabla \cdot(r\rho(\phi)\mathbf{u})\mathbf{u} \cdot \mathbf{u}+r\nabla p\cdot \mathbf{u}\right)drdz-\int_{\Omega}\nabla \cdot(r\eta(\phi)D(\mathbf{u}))\cdot\mathbf{u}drdz\nonumber\\
&-\frac{1}{2}\int_{\Omega}r\eta|D(\mathbf{u})|^{2}drdz +\int_{\Omega}\frac{1}{2}\nabla\cdot r\mathbf{J}(\mu)\mathbf{u}\cdot\mathbf{u}+r\mathbf{J}(\mu)\cdot\nabla\mathbf{u}\cdot\mathbf{u}drdz+\frac{\boldsymbol{K}\mathcal{K}}{\sqrt{\int_{0}^{t}-\mathcal{K}d\tau+\boldsymbol{G}}}\bigg), \label{229}\\
&\boldsymbol{K}_{t}=\frac{-\boldsymbol{R}\mathcal{K}}{2\sqrt{\int_{0}^{t}-\mathcal{K}d\tau+\boldsymbol{G}}},\label{230}\\
&\boldsymbol{T}_{t}=\alpha\int_{\Omega}\left(\frac{\partial (rv_{r})}{\partial r}+\frac{\partial (rv_{z})}{\partial z}\right)p drdz, \label{231}
\end{align}
where $\boldsymbol{H}(\phi)=(f(\phi)-s\phi)/\sqrt{\int_{\Omega}\left(rF(\phi)-\frac{rs}{2}\phi^{2}\right) drdz+B}$. $\alpha$ is a small enough positive number used to maintain the stability of the algorithm. Because $\boldsymbol{Q}_{t}=0$, $\boldsymbol{R}_{t}=0$ and $\boldsymbol{T}_{t}=0$, $\alpha$ does not change these equations.
Now we explain the modification made to the original system \eqref{27}-\eqref{210} to create a new system \eqref{224}-\eqref{231}. First, we incorporate the inner products of coupled nonlinear terms and certain non-zero boundary terms with specific functions into the ODE \eqref{228}, \eqref{229}, and \eqref{231}. Note all integral terms contained \eqref{228}, \eqref{229} and \eqref{231} are equal to zero, which means $\boldsymbol{Q}\equiv 1$, $\boldsymbol{R}\equiv 1$ and $\boldsymbol{T}\equiv 1$. Second, we add the term $r\frac{\mathcal{R}e}{2}\rho_{t}+\boldsymbol{R}\frac{\mathcal{R}e}{2}\nabla\cdot\left(r\rho\mathbf{u}\right)+\frac{\boldsymbol{R}}{2}\nabla\cdot r\mathbf{J}(\mu)$ to the momentum equation \eqref{29}. It is easy to know it is a zero from \eqref{27}. Third, we multiply certain terms by $\boldsymbol{Q}$ or $\boldsymbol{R}$. Since $\boldsymbol{Q}=\boldsymbol{R}=1$, these modifications will not alter the PDE system. Fourth, the ODE \eqref{227} is derived by taking the time derivative of $\boldsymbol{U}$. After integrating \eqref{227} with respect to time t and applying the initial condition, \eqref{28} is obtained. This implies \eqref{225} and \eqref{227} are equivalent to \eqref{28}. Therefore, the new PDE system \eqref{224}-\eqref{231} is equivalent to the original model \eqref{27}-\eqref{210}. 
Since the three equations \eqref{228}, \eqref{229}, and \eqref{231} are only differential equations with respect to time, the boundary conditions of the new system \eqref{224}-\eqref{231} remain the same as \eqref{212}. The transformed system \eqref{224}-\eqref{231} in the new variables $p, \mathbf{u}, \phi, \mu, \boldsymbol{U}, \boldsymbol{K},\boldsymbol{Q}, \boldsymbol{R}, \boldsymbol{T}$ forms a closed PDE system with the following initial conditions.
$$
\begin{aligned}
&(\phi,\mathbf{u},\mu)|_{t=0}=(1,\mathbf{0},0), \\
&\boldsymbol{U}|_{t=0}=\sqrt{\int_{\Omega}\big(rF(\phi^{0})-\frac{rs}{2}(\phi^{0})^{2}\big)drdz+B},\quad \boldsymbol{Q}|_{t=0}=1,\quad \boldsymbol{R}|_{t=0}=1,\quad \boldsymbol{T}|_{t=0}=1,\\
& \boldsymbol{K}|_{t=0}=\sqrt{G}.
\end{aligned}
$$
The modified system \eqref{224}-\eqref{231} also satisfies the law of energy dissipation.  It can be derived through a similar energy estimation process as \eqref{222}.
\begin{theorem}
The solutions of the modified system \eqref{224}-\eqref{231}, with boundary conditions \eqref{212}, satisfy
\begin{align}
&\frac{d}{dt}\mathbf{E}_{M}=-\mathcal{B}\boldsymbol{{L}_{d}}\int_{\Omega}r|\nabla \mu|^{2}drdz-\frac{\boldsymbol{R}}{2}\int_{\Omega}r\eta |\boldsymbol{D}(\mathbf{u})|^{2}drdz-2\int_{\Omega}\eta\frac{v_{r}^{2}}{r}drdz<0, \label{tcenergy}
\end{align}
where,
\begin{align}
\mathbf{E}_{M}=&\mathcal{R}e\frac{1}{2}\int_{\Omega}r|\sigma\mathbf{u}|^{2}drdz+\frac{\mathcal{B}\epsilon}{2}\int_{\Omega}r|\nabla \phi|^{2}drdz+\frac{\mathcal{B}s}{2\epsilon}\int_{\Omega}r|\phi|^{2}drdz+\frac{\mathcal{B}}{\epsilon}|\boldsymbol{U}|^{2}+\frac{1}{\alpha}\frac{\mathcal{B}|\boldsymbol{Q}|^{2}}{2}+\frac{1}{\alpha}\frac{|\boldsymbol{R}|^{2}}{2}+|\boldsymbol{K}|^{2}+\frac{1}{\alpha}\frac{|\boldsymbol{T}|^{2}}{2}. \label{233}
\end{align}
\end{theorem}
\noindent\textbf{Proof: }Taking the inner product of \eqref{224} with $r\mu$ in $L^{2}$, we have
\begin{align}
\int_{\Omega}r\phi_{t}\mu drdz+\boldsymbol{Q}\int_{\Omega}r\mathbf{u}\cdot \nabla\phi\mu drdz=-\mathcal{L}_{d}\int_{\Omega}r|\nabla \mu|^{2}drdz.
\end{align}
Taking the $L^{2}$ inner product of \eqref{225} with $r\phi_{t}$, we obtain
\begin{equation}
\int_{\Omega}r\phi_{t}\mu drdz=\frac{d}{dt}\left(\frac{\epsilon}{2}\int_{\Omega}r|\nabla \phi|^{2}drdz+\frac{s}{2\epsilon}\int_{\Omega}r|\phi|^{2}drdz\right)+\frac{\boldsymbol{U}}{\epsilon}\int_{\Omega}r\boldsymbol{H}\phi_{t}drdz.
\end{equation}
Combining the above two equations, we have
\begin{align}
\frac{d}{dt}\left(\frac{\epsilon}{2}\int_{\Omega}r|\nabla \phi|^{2}drdz+\frac{s}{2\epsilon}\int_{\Omega}r|\phi|^{2}drdz\right)=&-\boldsymbol{Q}\int_{\Omega}r\mathbf{u}\cdot \nabla\phi\mu drdz-\frac{\boldsymbol{U}}{\epsilon}\int_{\Omega}r\boldsymbol{H}\phi_{t}drdz-\mathcal{L}_{d}\int_{\Omega}r|\nabla \mu|^{2}drdz.
\end{align}
Taking the $L^{2}$ inner product of \eqref{227} with $2\frac{\boldsymbol{U}}{\epsilon}$, we obtain 
\begin{equation}
\frac{1}{\epsilon}\frac{d}{dt}|\boldsymbol{U}|^{2}=\frac{\boldsymbol{U}}{\epsilon}\int_{\Omega}r\boldsymbol{H}\phi_{t}drdz.
\end{equation}
Combining the above two equations and multiply it with $\mathcal{B}$, then we have
\begin{align}
&\frac{d}{dt}\left(\frac{\mathcal{B}\epsilon}{2}\int_{\Omega}r|\nabla \phi|^{2}drdz+\frac{\mathcal{B}s}{2\epsilon}\int_{\Omega}r|\phi|^{2}drdz\right)+\frac{\mathcal{B}}{\epsilon}\frac{d}{dt}|\boldsymbol{U}|^{2}=-\mathcal{B}\boldsymbol{Q}\int_{\Omega}r\mathbf{u}\cdot\nabla\phi\mu drdz-\mathcal{B}\mathcal{L}_{d}\int_{\Omega}r|\nabla \mu|^{2}drdz.
\end{align}
Multiplying \eqref{228} with $\frac{\mathcal{B}}{\alpha}\boldsymbol{Q}$, we obtain
\begin{align}
\frac{\mathcal{B}}{\alpha}\frac{d}{dt}(\frac{|\boldsymbol{Q}|^{2}}{2})=\boldsymbol{Q}\mathcal{B}\int_{\Omega}\bigg(r\mathbf{u}\cdot\left(\frac{\partial \phi}{\partial z},\frac{\partial \phi}{\partial r}\right) \mu-r\mathbf{u}\cdot\left(\frac{\partial \phi}{\partial z},\frac{\partial \phi}{\partial r}\right) \mu\bigg) drdz.
\end{align}
Taking the  $L^{2}$ inner product of \eqref{226} with $r\mathbf{u}$ and using integration by parts and the divergence free condition \eqref{232}, we obtain
\begin{equation}
    \begin{split}
&\mathcal{R}e\frac{1}{2}\frac{d}{dt}\int_{\Omega}r|\sigma\mathbf{u}|^{2}drdz+\int_{\Omega}\left(\frac{\mathcal{R}e\boldsymbol{R}}{2}\nabla\cdot(r\rho\mathbf{u})\mathbf{u}\cdot\mathbf{u}+\mathcal{R}e\boldsymbol{R}r\rho(\mathbf{u}\cdot\nabla)\mathbf{u}\cdot\mathbf{u}\right)drdz\\
&+(1-\mathbf{R})\int_{\Omega}\nabla\cdot\left( r\eta\boldsymbol{D}(\mathbf{u})\right)\cdot \mathbf{u} drdz+\boldsymbol{R}\int_{\Omega}r\nabla p\cdot\mathbf{u}drdz-\boldsymbol{Q}\mathcal{B}\int_{\Omega}r\mu\nabla\phi\cdot\mathbf{u}drdz.\\
&+2\int_{\Omega}\eta\frac{v_{r}^{2}}{r}drdz-\int_{\Omega}\nabla\cdot\left( r\eta\boldsymbol{D}(\mathbf{u})\right)\cdot \mathbf{u} drdz+\boldsymbol{R}\int_{\Omega}\nabla\cdot r\mathbf{J}(\mu)\mathbf{u}\cdot\mathbf{u}+r\mathbf{J}(\mu)\cdot\nabla\mathbf{u}\cdot\mathbf{u}drdz=0.
    \end{split}
\end{equation}
Multiplying \eqref{229} with $\frac{1}{\alpha}\boldsymbol{R}$, we obtain
\begin{equation}
    \begin{split}
\frac{1}{\alpha}\frac{d}{dt}\left(\frac{|\boldsymbol{R}|^{2}}{2}\right)=&\boldsymbol{R}\int_{\Omega}\big(\mathcal{R}e r\rho(\phi)(\mathbf{u} \cdot \nabla) \mathbf{u} \cdot \mathbf{u}+\frac{\mathcal{R}e}{2}\nabla \cdot(r\rho(\phi)\mathbf{u})\mathbf{u} \cdot \mathbf{u}+r\nabla p\cdot \mathbf{u}\big)drdz.\\
&-\boldsymbol{R}\int_{\Omega}\nabla \cdot(r\eta(\phi)D(\mathbf{u}))\cdot\mathbf{u}drdz-\frac{\boldsymbol{R}}{2}\int_{\Omega}r\eta|D(\mathbf{u})|^{2}drdz+\\
&\boldsymbol{R}\int_{\Omega}\nabla\cdot r\mathbf{J}(\mu)\mathbf{u}\cdot\mathbf{u}+r\mathbf{J}(\mu)\cdot\nabla\mathbf{u}\cdot\mathbf{u}drdz+\frac{\boldsymbol{RK}\mathcal{K}}{\sqrt{\int_{0}^{t}-\mathcal{K} d\tau+\boldsymbol{G}}}.
    \end{split}
\end{equation}
Multiplying \eqref{230} with 2$\boldsymbol{K}$, we obtain
\begin{align}
&\frac{d}{dt}(|\boldsymbol{K}|^{2})=-\frac{\boldsymbol{RK}\mathcal{K}}{\sqrt{\int_{0}^{t}-\mathcal{K}d\tau+\boldsymbol{G}}}.
\end{align}
Multiplying \eqref{231} with $\frac{1}{\alpha}\boldsymbol{T}$, we obtain
\begin{align}
\frac{1}{\alpha}\frac{d}{dt}(\frac{|\boldsymbol{T}|^{2}}{2})=\boldsymbol{T}\int_{\Omega}\left(\frac{\partial (rv_{r})}{\partial r}+\frac{\partial (rv_{z})}{\partial z}\right)p drdz.
\end{align}
By combining the above equations, we obtain the law of energy dissipation of the modified system \eqref{224}-\eqref{231} as follows
\begin{equation}
    \begin{split}
&\frac{d}{dt}\bigg(\mathcal{R}e\frac{1}{2}\int_{\Omega}r|\sigma\mathbf{u}|^{2}drdz+\frac{\mathcal{B}\epsilon}{2}\int_{\Omega}r|\nabla \phi|^{2}drdz+\frac{\mathcal{B}s}{2\epsilon}\int_{\Omega}r|\phi|^{2}drdz+\frac{\mathcal{B}}{\epsilon}|\boldsymbol{U}|^{2}+\frac{1}{\alpha}\frac{\mathcal{B}|\boldsymbol{Q}|^{2}}{2}\\
&+\frac{1}{\alpha}\frac{|\boldsymbol{R}|^{2}}{2}+|\boldsymbol{K}|^{2}+\frac{1}{\alpha}\frac{|\boldsymbol{T}|^{2}}{2}\bigg)=-\mathcal{B}\boldsymbol{{L}_{d}}\int_{\Omega}r|\nabla \mu|^{2}drdz-\frac{\boldsymbol{R}}{2}\int_{\Omega}r\eta |\boldsymbol{D}(\mathbf{u})|^{2}drdz-2\int_{\Omega}\eta\frac{v_{r}^{2}}{r}drdz<0.
    \end{split}
\end{equation}

It is worth noting that the modified energy $E_{M}$ and the original energy $E_{O}$ are equivalent in the sense of a constant difference. \begin{align}
 E_{M} - E_{O}=\frac{1}{\alpha}\frac{\mathcal{B}|\boldsymbol{Q}|^{2}}{2}+\frac{1}{\alpha}\frac{|\boldsymbol{R}|^{2}}{2}+\boldsymbol{G}+\frac{\mathcal{B}}{\epsilon}B+\frac{1}{\alpha}\frac{|\boldsymbol{T}|^{2}}{2}\label{dE}
\end{align}
can be deduced by $\boldsymbol{U}=\sqrt{\int_{\Omega}\left(rF(\phi)-\frac{rs}{2}\phi^{2}\right) drdz+B}$, $\frac{\mathcal{B}s}{2\epsilon}\int_{\Omega}r|\phi|^{2}drdz+\frac{\mathcal{B}}{\epsilon}|\boldsymbol{U}|^{2}=\frac{\mathcal{B}}{\epsilon}\int_{\Omega}rF(\phi)drdz+\frac{\mathcal{B}}{\epsilon}B$ and $\boldsymbol{K}=\sqrt{\int_{0}^{t}-\mathcal{K}d\tau+\boldsymbol{G}}$. Now we can see that $E_{M} - E_{O}$ is a constant.

\section{Numerical scheme and its solvability}
\subsection{Numerical scheme and its energy law}
In this subsection, we present a time-marching scheme to solve the system \eqref{224}-\eqref{231} which is an equivalent system of \eqref{27}-\eqref{210}. We denote $\delta_{t}>0$ as a time step size and $t^{n}=n_{t} (0\le n\le N)$ with $T = N_{t}$. Let $\psi^{n}$ represent the numerical approximation to the function  $\psi(\cdot,t)|_{t=t^{n}}$.
$$
\begin{aligned}
\mathcal{K}^{n}=&\frac{\mathcal{R}e}{2}\int_{\partial\Omega_{L}}r\rho^{n+1}({v}^{n}_{z})^{3}d \boldsymbol{s}-\frac{\mathcal{R}e}{2}\int_{\partial\Omega_{R}}r\rho^{n+1} ({v}^{n}_{z})^{3}d \boldsymbol{s}+\int_{\partial\Omega_{L}}r(2p^{n}-p^{n-1})v^{n}_{z}d s\\
&-2\int_{\partial\Omega_{L}}r\eta^{n+1} v^{n}_{z}\frac{\partial v^{n}_{z}}{\partial z}ds-\int_{\partial\Omega_{R}}r(2p^{n}-p^{n-1})v^{n}_{z}ds+\mathcal{L}_{d}\mathcal{R}e\frac{1-\lambda_{\rho}}{4}\int_{\partial\Omega_{L}}\frac{\partial\mu^{n}}{\partial z}(v_{z}^{n})^{2}ds.
\end{aligned}
$$
Since the boundary condition of $\mathbf{u}$ is non-homogeneous Dirichlet boundary condition, if we use the penalty method \cite{ref29} for the Navier-Stokes equation in our problem, it is difficult to prove the energy law of the numerical scheme. To address this challenge, we introduce a nonlocal variable $\boldsymbol{T}$ and present a modified penalty method that overcomes the difficulty. The Strang operator splitting method \cite{ref26,ref27} is used to decouple the momentum equation and the phase-field equation. SAV method\cite{ref24,ref25} is utilized to linearize the nonlinear term $f(\phi)$ in the Cahn-Hilliard equation. A scheme to solve \eqref{224}-\eqref{231} is constructed as follows: Given $\phi^{n}$, $\mathbf{u}^{n}$, $\mu^{n}$, $\boldsymbol{U}^{n}$, $\boldsymbol{Q}^{n}$, $\boldsymbol{R}^{n}$, $\boldsymbol{K}^{n}$, $p^{n}$, we compute $\phi^{n+1}$, $\mathbf{u}^{n+1}$, $\mu^{n+1}$, $\boldsymbol{U}^{n+1}$, $\boldsymbol{Q}^{n+1}$, $\boldsymbol{R}^{n+1}$, $\boldsymbol{K}^{n+1}$, $p^{n+1}$ by the following three steps.\\
$\boldsymbol{STEP 1:}$
\begin{align}
&r\frac{\phi^{n+1}-\phi^{n}}{\delta_{t}}+\boldsymbol{Q}^{n+1}r\mathbf{u}^{n}\cdot\nabla\phi^{n}=\mathcal{L}_{d}\nabla\cdot(r\nabla\mu^{n+1}), \label{31}\\
&r\mu^{n+1}=-\epsilon \nabla\cdot(r\nabla\phi^{n+1})+\frac{rs}{\epsilon}\phi^{n+1}+\frac{r}{\epsilon}\boldsymbol{H}^{n}\boldsymbol{U}^{n+1}, \label{32}\\
&\mathcal{R}e\rho^{n}\frac{\tilde{\mathbf{u}}^{n+1}-\mathbf{u}^{n}}{\delta_{t}}-\boldsymbol{Q}^{n+1}\mathcal{B}\mu^{n}\nabla\phi^{n}=0,\label{33}\\
&\boldsymbol{U}^{n+1}-\boldsymbol{U}^{n}=\frac{1}{2}\int_{\Omega}r\boldsymbol{H}^{n}(\phi^{n+1}-\phi^{n})drdz, \label{34}\\
&\frac{\boldsymbol{Q}^{n+1}-\boldsymbol{Q}^{n}}{\delta_{t}}=\alpha\int_{\Omega}(r\mathbf{u}^{n}\cdot\nabla\phi^{n} \mu^{n+1}-r\tilde{\mathbf{u}}^{n+1}\cdot\nabla\phi^{n} \mu^{n}) drdz. \label{35}
\end{align}
The boundary conditions read as
\begin{equation}
    \begin{split}
&\frac{\partial\mu^{n+1}}{\partial\mathbf{n}}(z, r, t)=0,\quad\frac{\partial\phi^{n+1}}{\partial\mathbf{n}}(z, r, t)=0,\quad \left\{(z, r)|(z, r)\in\Gamma_{1}\cup\Gamma_{4}\cup\Gamma_{5}\right\},\\
&\phi^{n+1}(z, r, t)=1, \quad\mu^{n+1}(z, r, t)=0, \quad\left\{(z, r)|(z, r)\in\Gamma_{2}\right\},\\
&\phi^{n+1}(z, r, t)=-1,\quad \mu^{n+1}(z, r, t)=0,\quad\left\{(z, r)|(z, r)\in\Gamma_{3}\right\}.\\
    \end{split}
\end{equation}
$\boldsymbol{STEP 2:}$
\begin{align}
&\mathcal{R}e\big[\frac{r}{2}\frac{\rho^{n+1}-\rho^{n}}{\delta_{t}}\mathbf{u}^{n+1}+r\rho^{n}\frac{\mathbf{u}^{n+1}-\tilde{\mathbf{u}}^{n+1}}{\delta_{t}}+\boldsymbol{R}^{n+1}r\rho^{n}(\mathbf{u}^{n}\cdot\nabla) \mathbf{u}^{n}+\frac{\boldsymbol{R}^{n+1}}{2}\nabla\cdot(r\rho^{n+1}\mathbf{u}^{n})\mathbf{u}^{n}\big]\nonumber\\ 
&-\nabla\cdot\left(r\eta^{n+1}D(\mathbf{u}^{n+1})\right)+\left(0, 2\eta^{n+1}\frac{v_{r}^{n+1}}{r}\right) +r\boldsymbol{R}^{n+1}\nabla(2p^{n}-p^{n-1})\label{36}\\ \nonumber
&+\nabla \cdot\left(r(\sqrt{\eta^{n}\eta^{n+1}}-\boldsymbol{R}^{n+1}\eta^{n+1})D(\mathbf{u}^{n})\right)+\frac{\boldsymbol{R}^{n+1}}{2}\nabla\cdot r\mathbf{J}(\mu^{n})\mathbf{u}^{n}+r\boldsymbol{R}^{n+1}\mathbf{J}(\mu^{n})\cdot\nabla\mathbf{u}^{n}=0, \\
&\frac{\boldsymbol{R}^{n+1}-\boldsymbol{R}^{n}}{\delta_{t}}=\alpha\bigg(\int_{\Omega}\left(\mathcal{R}er\rho^{n+1}(\mathbf{u}^{n}\cdot\nabla) \mathbf{u}^{n} \cdot \mathbf{u}^{n+1}+\frac{\mathcal{R}e}{2}\nabla\cdot(r\rho^{n+1}\mathbf{u}^{n})\mathbf{u}^{n} \cdot \mathbf{u}^{n+1}\right) drdz \nonumber\\
&+\int_{\Omega}r\nabla(2p^{n}-p^{n-1})\cdot \mathbf{u}^{n+1}drdz-\int_{\Omega}\nabla\cdot\left(r\eta^{n+1}D(\mathbf{u}^{n})\right)\cdot\mathbf{u}^{n+1}drdz-\frac{\boldsymbol{R}^{n+1}}{2}\int_{\Omega}r\eta^{n+1}|D(\mathbf{u}^{n})|^{2}drdz\label{37}\\ \nonumber
&+\int_{\Omega}\frac{1}{2}\nabla\cdot r\mathbf{J}(\mu^{n})\mathbf{u}^{n}\cdot\mathbf{u}^{n+1}+r\mathbf{J}(\mu^{n})\cdot\nabla\mathbf{u}^{n}\cdot\mathbf{u}^{n+1}drdz+\frac{\boldsymbol{K}^{n+1}\mathcal{K}^{n}}{\sqrt{\int_{0}^{t_{n}}-\mathcal{K}d\tau+\boldsymbol{G}}}\bigg),\\
&\frac{\boldsymbol{K}^{n+1}-\boldsymbol{K}^{n}}{\delta_{t}}=\frac{-\boldsymbol{R}^{n+1}\mathcal{K}^{n}}{2\sqrt{\int_{0}^{t_{n}}-\mathcal{K}d\tau+\boldsymbol{G}}}. \label{38}
\end{align}
The boundary conditions read as
\begin{equation}
    \begin{split}
&\mathbf{u}^{n+1}(z, r, t)=0,\quad \left\{(z, r)|(z, r)\in\Gamma_{1}\right\},\\
&v_r^{n+1}(z, r, t)=0, \quad v_z^{n+1}(z, r, t)=2 \frac{Q_r}{a^2}\left[\frac{1-\left(\frac{r}{a}\right)^2+\frac{1-\left(1 /a\right)^2}{\ln \left(a\right)} \ln \left(\frac{r}{a}\right)}{1-\left(\frac{1}{a}\right)^4-\frac{\left(1-\left(1 / a\right)^2\right)^2}{\ln \left(a\right)}}\right],\quad\left\{(z, r)|(z, r)\in\Gamma_{2}\right\},\\
&v_r^{n+1}(z, r, t)=0, \quad v_z^{n+1}(z, r, t)=2\left[1-r^2\right],\left\{(z, r)|(z, r)\in\Gamma_{3}\right\},\\
&v_r^{n+1}(z, r, t)=0, \quad \frac{\partial v_z^{n+1}}{\partial\mathbf{n}}(z, r, t)=0,\quad\left\{(z, r)|(z, r)\in\Gamma_{4}\cup\Gamma_{5}\right\}.\\
    \end{split}
\end{equation}
$\boldsymbol{STEP 3:}$
\begin{align}
&\frac{\boldsymbol{T}^{n+1}-\boldsymbol{T}^{n}}{\delta_{t}}=\alpha\int_{\Omega}\nabla\cdot(r\mathbf{u}^{n+1})p^{n+1}drdz, \label{39}\\
&\nabla\cdot \left(r\nabla(p^{n+1}-p^{n})\right)=\boldsymbol{T}^{n+1}\frac{\chi\mathcal{R}e}{\delta_{t}}\nabla\cdot(r\mathbf{u}^{n+1}). \label{310}
\end{align}
The boundary conditions read as
\begin{equation}
    \begin{split}
&\frac{\partial p^{n+1}}{\partial\mathbf{n}}(z, r, t)=0,\quad\left\{(z, r)|(z, r)\in\Gamma_{1}\cup\Gamma_{2}\cup\Gamma_{3}\cup\Gamma_{4}\right\},\\
&p^{n+1}=0, \quad\left\{(z, r)|(z, r)\in\Gamma_{5}\right\}.\label{pb}
    \end{split}
\end{equation}
Other notations used in the scheme read as
$$
\left\{\begin{array}{l}
H^n=H\left(\phi^n\right),\\
\chi=\frac{1}{2} \min \left(\rho_1, \rho_2\right), \\
\hat{\phi}= \begin{cases}\phi, & |\phi|<1, \\
\operatorname{sign}(\phi), & |\phi|>1,\end{cases} \\
\rho^{n+1}=\frac{\lambda_{\rho}-1}{2} \hat{\phi}^{n+1}+\frac{\lambda_{\rho}+1}{2},\quad\eta^{n+1}=\frac{\lambda_{\eta}-1}{2} \hat{\phi}^{n+1}+\frac{\lambda_{\eta}+1}{2} .
\end{array}\right.
$$

In the last section, we have mentioned that the aim of adding the term $\frac{(1-\mathbf{R})}{r}\nabla \cdot(r\eta(\phi)D(\mathbf{u}))$ is to enhance the stability of our numerical scheme. Because the term $\nabla\cdot(r\eta D(\mathbf{u}))$ will produce the non-zero boundary terms during the process of proving the energy law, we need to use the nonlocal variable $\boldsymbol{R}$ and $\boldsymbol{K}$ to eliminate these terms. But the disadvantage of this approach is that we need to process explicitly the term $\nabla\cdot(r\eta D(\mathbf{u}))$. It will make the numerical scheme unstable, that is the value of the nonlocal variable $\boldsymbol{Q}$, $\boldsymbol{R}$, $\boldsymbol{T}$ will deviate from 1 during calculation, so we introduce the term $\frac{(1-\mathbf{R})}{r}\nabla \cdot(r\eta(\phi)D(\mathbf{u}))$ and process this term explicitly and the term $\nabla\cdot(r\eta D(\mathbf{u}))$ implicitly to overcome this problem, and then we can prove the numerical scheme is unconditional energy stable.\par
The scheme is linear, and all nonlinear terms are discretized using a combination of explicit and implicit methods, even though there are some terms involved in $\boldsymbol{STEP 2}$, such as $\nabla\cdot(r\eta^{n+1}D(\mathbf{u}^{n+1}))$, $\nabla \cdot\left(r(\sqrt{\eta^{n}\eta^{n+1}}-\boldsymbol{R}^{n+1}\eta^{n+1})D(\mathbf{u}^{n})\right)$ and $r\rho^{n+1}\frac{\mathbf{u}^{n+1}-\mathbf{u}^{n}}{\delta_{t}}$ in \eqref{36}. Considering that $\eta^{n+1}$ and $\rho^{n+1}$ in $\boldsymbol{STEP 2}$ are the known terms since $\phi^{n+1}$ is obtained in $\boldsymbol{STEP 1}$, the scheme in $\boldsymbol{STEP 2}$ is still linear . In $\boldsymbol{STEP 2}$, we use the first-order operator Strang splitting method\cite{ref22} to split the surface tension term $\boldsymbol{Q}^{n+1}\phi^{n}\nabla\mu^{n}$ from the momentum equation to obtain the intermediate velocity $\tilde{\mathbf{u}}^{n+1}$. It is evident that $\boldsymbol{Q}^{n+1}$, $\boldsymbol{R}^{n+1}$, $\boldsymbol{T}^{n+1}$ will not retain their exact value 1 in the calculation, because $\boldsymbol{Q}^{n+1}$, $\boldsymbol{R}^{n+1}$, $\boldsymbol{T}^{n+1}$ are only numerical approximations to $\boldsymbol{Q}$, $\boldsymbol{R}$, $\boldsymbol{T}$. Therefore, we multiply the right-hand side of ODEs by $\alpha$ which is a small enough number to maintain the stability of the numerical scheme. The accuracy test also illustrates this fact.

The following theorem shows the unconditional energy stability of the numerical scheme.
\begin{theorem}
The solutions of the time-discrete scheme \eqref{31}-\eqref{pb} satisfy
\begin{align}
\mathbf{E}^{n+1}_{M}-\mathbf{E}^{n}_{M}<&-\delta_t \mathcal{B} \mathcal{L}_d \mathcal{R}e \int_{\Omega}\left|\nabla r \mu^{n+1}\right|^2 drdz-\frac{\delta_t |R_{n+1}|^{2}}{2} \int_{\Omega} r \eta^{n+1}\left|D\left(\mathbf{u}^{n+1}\right)\right|^2 drdz-\delta_t\int_{\Omega} 2\eta^{n+1}\frac{(v_{r}^{n+1})^{2}}{r}drdz, \label{tdenergy}
\end{align}
where,
\begin{equation}
    \begin{split}
\mathbf{E}^{n+1}_{M}=& \frac{\mathcal{R}e}{2} \int_{\Omega} r \rho^{n+1}\left|\mathbf{u}^{n+1}\right|^2 drdz+\frac{\mathcal{B}}{2} \epsilon \int_{\Omega}r\left|\nabla \phi^{n+1}\right|^2drdz+\frac{\mathcal{B}\boldsymbol{s}}{2 \epsilon} \int_{\Omega}r\left|\phi^{n+1}\right|^2drdz \\ 
&+2\delta_{t}\int_{\partial\Omega_{L}}r\eta^{n+1} v^{n+1}_{z}\frac{\partial v^{n+1}_{z}}{\partial z}ds+\frac{\delta_{t}}{4}\int_{\Omega}r\eta^{n+1}|D(\mathbf{u}^{n+1})|^{2}drdz+\frac{1}{\alpha}\frac{\mathcal{B}}{2}\left|\boldsymbol{Q}^{n+1}\right|^2\\ 
&+\frac{\mathcal{B}}{\epsilon}\left|U^{n+1}\right|^2+\frac{1}{\alpha}\frac{\left|\boldsymbol{R}^{n+1}\right|^2}{2}+\frac{1}{\alpha}\frac{\left|\boldsymbol{T}^{n+1}\right|^2}{2}+\left|\boldsymbol{K}^{n+1}\right|^2+\frac{\delta_{t}^{2}}{2\chi\mathcal{R}e} \int_{\Omega} r\left|\nabla p^{n+1}\right|^2drdz.
    \end{split}\label{EM}
\end{equation}
\end{theorem}
\noindent \textbf{Proof:}
We take the $L^{2}$ inner product of \eqref{31} with $\mu^{n+1} _t$ and use integration by parts to get
\begin{equation}
\int_{\Omega} r \mu^{n+1}\left(\phi^{n+1}-\phi^n\right) drdz+\boldsymbol{Q}^{n+1} _t \int_{\Omega} r (\mathbf{u}^n \cdot \nabla \phi^n) \mu^{n+1} drdz=-\delta_t \mathcal{L}_d \int_{\Omega} r\left|\nabla \mu^{n+1}\right|^2 drdz.
\end{equation}
Taking the $L^{2}$ inner product of \eqref{32} with  $\left(\phi^{n+1}-\phi^n\right)$ and using integration by parts, we
get
\begin{equation}
    \begin{split}
& \int_{\Omega} r\left(\phi^{n+1}-\phi^n\right) \mu^{n+1} drdz =\frac{1}{2} \epsilon \int_{\Omega} \left(r\left|\nabla \phi^{n+1}\right|^2-r\left|\nabla \phi^n\right|^2+r\left|\nabla \phi^{n+1}-\nabla \phi^n\right|^2 \right) drdz \\
& +\frac{s}{2 \epsilon} \int_{\Omega}\left( r\left|\phi^{n+1}\right|^2-r\left|\phi^n\right|^2+r\left|\phi^{n+1}-\phi^n\right|^2\right) drdz+\frac{\boldsymbol{U}^{n+1}}{\epsilon}  \int_{\Omega} r \boldsymbol{H}^n\left(\phi^{n+1}-\phi^n\right) drdz.
    \end{split}
\end{equation}
Combining the above two equations and multiplying it with $\mathcal{B}$, we obtain
\begin{equation}
    \begin{split}
&\frac{\mathcal{B}}{2} \epsilon \int_{\Omega} \left(r\left|\nabla \phi^{n+1}\right|^2-r\left|\nabla \phi^n\right|^2+r\left|\nabla \phi^{n+1}-\nabla \phi^n\right|^2 \right) drdz \\
& +\frac{\mathcal{B}s}{2 \epsilon} \int_{\Omega}\left( r\left|\phi^{n+1}\right|^2-r\left|\phi^n\right|^2+r\left|\phi^{n+1}-\phi^n\right|^2\right) drdz+\mathcal{B}\frac{\boldsymbol{U}^{n+1}}{\epsilon} \int_{\Omega} r \boldsymbol{H}^n\left(\phi^{n+1}-\phi^n\right) drdz\\
&+\mathcal{B}\boldsymbol{Q}^{n+1} _t \int_{\Omega} r (\mathbf{u}^n \cdot \nabla \phi^n) \mu^{n+1} drdz=-\delta_t\mathcal{B} \mathcal{L}_d \int_{\Omega} r\left|\nabla \mu^{n+1}\right|^2 drdz.
    \end{split}
\end{equation}
Taking the $L^{2}$ inner product of \eqref{33} with $r\tilde{\mathbf{u}}^{n+1} _t$, we obtain
\begin{equation}
\frac{\mathcal{R}e}{2} \int_{\Omega} r \rho^n\left(\left|\tilde{\mathbf{u}}^{n+1}\right|^2-|\mathbf{u}^{n}|^2+\left|\tilde{\mathbf{u}}^{n+1}-\mathbf{u}^{n}\right|^2\right) drdz-\delta_t \boldsymbol{Q}^{n+1}\mathcal{B} \int_{\Omega} r \mu^n \nabla \phi^{n}\cdot \tilde{\mathbf{u}}^{n+1} drdz=0.
\end{equation}
We multiply \eqref{34} with $2 \frac{\mathcal{B}}{\epsilon} \boldsymbol{U}^{n+1}$ to get
\begin{equation}
\frac{\mathcal{B}}{\epsilon}\left|\boldsymbol{U}^{n+1}\right|^2-\frac{\mathcal{B}}{\epsilon}\left|\boldsymbol{U}^n\right|^2+\frac{\mathcal{B}}{\epsilon}\left|\boldsymbol{U}^{n+1}-\boldsymbol{U}^n\right|^2=\frac{\mathcal{B}}{\epsilon} \boldsymbol{U}^{n+1} \int_{\Omega} r \boldsymbol{H}^n\left(\phi^{n+1}-\phi^n\right) drdz.
\end{equation}
We multiply \eqref{35} with $\frac{\mathcal{B}}{\alpha}\boldsymbol{Q}^{n+1} _t$ to get
\begin{equation}
\frac{\mathcal{B}}{\alpha}\frac{\left|\boldsymbol{Q}^{n+1}\right|^2}{2}-\frac{\mathcal{B}}{\alpha}\frac{\left|\boldsymbol{Q}^n\right|^2}{2}+\frac{\mathcal{B}}{\alpha}\frac{\left|\boldsymbol{Q}^{n+1}-\boldsymbol{Q}^n\right|^2}{2}=\mathcal{B}\boldsymbol{Q}^{n+1} _t \int_{\Omega} r \mathbf{u}^n \cdot \nabla \phi^n \mu^{n+1}-r \tilde{\mathbf{u}}^{n+1} \cdot \nabla \phi^n \mu^n drdz.
\end{equation}
Taking the inner product of \eqref{36} with $\mathbf{u}^{n+1} _t$ in the $L^{2}$ space and using the following identity, we have
\begin{equation}
    \begin{split}
& \mathcal{R}e\left[\int_{\Omega}\left(\frac{r}{2}\left(\rho^{n+1}-\rho^n\right)\left|\mathbf{u}^{n+1}\right|^2+\frac{r}{2} \rho^n\left(\left|\mathbf{u}^{n+1}\right|^2-\left|\tilde{\mathbf{u}}^{n+1}\right|^2+\left|\mathbf{u}^{n+1}-\tilde{\mathbf{u}}^{n+1}\right|^2\right)\right) drdz+\right. \\
& \left.\delta_t \boldsymbol{R}^{n+1} \int_{\Omega}\big( \frac{1}{2} \nabla \cdot\left(r \rho^n \mathbf{u}^n\right) \mathbf{u}^n \cdot \mathbf{u}^{n+1}+r \rho^n\left(\mathbf{u}^n \cdot \nabla\right) \mathbf{u}^n \cdot \mathbf{u}^{n+1}\big) drdz\right]+\delta_t \boldsymbol{R}^{n+1} \int_{\Omega} r \nabla\big(2 p^n-p^{n-1}\big) \cdot \mathbf{u}^{n+1} drdz\\
&-\delta_t \int_{\Omega} \nabla \cdot\left(r \eta^{n+1} \boldsymbol{D}\left(\mathbf{u}^{n+1}\right)-r\sqrt{\eta^{n}\eta^{n+1}}\boldsymbol{D}\left(\mathbf{u}^{n}\right)\right) \cdot\mathbf{u}^{n+1}drdz-\delta_t\boldsymbol{R}^{n+1} \int_{\Omega} \nabla \cdot\left(r \eta^{n+1} \boldsymbol{D}\left(\mathbf{u}^{n}\right)\right) \cdot\mathbf{u}^{n+1}drdz\\
&+2\delta_{t}\int_{\Omega}\eta^{n+1}\frac{(v_{r}^{n+1})^{2}}{r}drdz+\delta_{t}\boldsymbol{R}^{n+1}\int_{\Omega}\frac{1}{2}\nabla\cdot r\mathbf{J}(\mu^{n})\mathbf{u}^{n}\cdot\mathbf{u}^{n+1}+r\mathbf{J}(\mu^{n})\cdot\nabla\mathbf{u}^{n}\cdot\mathbf{u}^{n+1}drdz=0.
    \end{split}
\end{equation}
The term $\int_{\Omega} \nabla \cdot\left(r \eta^{n+1} \boldsymbol{D}\left(\mathbf{u}^{n+1}\right)-r\sqrt{\eta^{n}\eta^{n+1}}\boldsymbol{D}\left(\mathbf{u}^{n}\right)\right) \cdot\mathbf{u}^{n+1}drdz$ is equal to 
$$
\begin{aligned}
&\int_{\Omega} \nabla \cdot\left(r \eta^{n+1} \boldsymbol{D}\left(\mathbf{u}^{n+1}\right)-r\sqrt{\eta^{n}\eta^{n+1}}\boldsymbol{D}\left(\mathbf{u}^{n}\right)\right) \cdot\mathbf{u}^{n+1}drdz\\
=&\int_{\partial\Omega}\mathbf{n}\cdot\left(r \eta^{n+1} \boldsymbol{D}\left(\mathbf{u}^{n+1}\right)-r\sqrt{\eta^{n}\eta^{n+1}}\boldsymbol{D}\left(\mathbf{u}^{n}\right)\right) \cdot\mathbf{u}^{n+1}ds-\int_{\Omega}\left(r \eta^{n+1} \boldsymbol{D}\left(\mathbf{u}^{n+1}\right)-r\sqrt{\eta^{n}\eta^{n+1}}\boldsymbol{D}\left(\mathbf{u}^{n}\right)\right):\nabla\mathbf{u}^{n+1}drdz\\
=&-\int_{\partial\Omega_{L}}2r\eta^{n+1}\left(v_{z}^{n+1}\frac{\partial \left(v_{z}^{n+1}-v_{z}^{n}\right)}{\partial z}+\left(v_{z}^{n+1}- v_{z}^{n}\right)\frac{\partial v_{z}^{n+1}}{\partial z}\right)ds\\
&-\int_{\Omega}\left(r \eta^{n+1} \boldsymbol{D}\left(\mathbf{u}^{n+1}\right)-r\sqrt{\eta^{n}\eta^{n+1}}\boldsymbol{D}\left(\mathbf{u}^{n}\right)\right):\nabla\mathbf{u}^{n+1}drdz\\
=&-\int_{\partial\Omega_{L}}2r\eta^{n+1}\frac{\partial \left(v_{z}^{n+1}(v_{z}^{n+1}-v_{z}^{n})\right)}{\partial z}ds-\int_{\Omega}\left(r \sqrt{\eta^{n+1}} \boldsymbol{D}\left(\mathbf{u}^{n+1}\right)-r\sqrt{\eta^{n}}\boldsymbol{D}\left(\mathbf{u}^{n}\right)\right):\sqrt{\eta^{n+1}}\nabla\mathbf{u}^{n+1}drdz\\
=&-\int_{\partial\Omega_{L}}r\eta^{n+1}\frac{\partial \left(|v_{z}^{n+1}|^{2}-|v_{z}^{n}|^{2}+|v_{z}^{n+1}-v_{z}^{n}|^{2}\right)}{\partial z}ds-\frac{1}{4}\int_{\Omega}\bigg(r |\sqrt{\eta^{n+1}} \boldsymbol{D}\left(\mathbf{u}^{n+1}\right)|^{2}-r|\sqrt{\eta^{n}}\boldsymbol{D}\left(\mathbf{u}^{n}\right)|^{2}\bigg)drdz\\
=&-\left(\int_{\partial\Omega_{L}}2r\eta^{n+1} v_{z}^{n+1}\frac{\partial v_{z}^{n+1}}{\partial z}ds-\int_{\partial\Omega_{L}}2r\eta^{n+1} v_{z}^{n}\frac{\partial v_{z}^{n}}{\partial z}ds\right)-\frac{1}{4}\int_{\Omega}\bigg(r |\sqrt{\eta^{n+1}} \boldsymbol{D}\left(\mathbf{u}^{n+1}\right)|^{2}-r|\sqrt{\eta^{n}}\boldsymbol{D}\left(\mathbf{u}^{n}\right)|^{2}\bigg)drdz\\
=&-\left(\int_{\partial\Omega_{L}}2r\eta^{n+1} v_{z}^{n+1}\frac{\partial v_{z}^{n+1}}{\partial z}ds-\int_{\partial\Omega_{L}}2r\eta^{n+1} v_{z}^{n}\frac{\partial v_{z}^{n}}{\partial z}ds\right)-\frac{1}{4}\int_{\Omega}\bigg(r |\sqrt{\eta^{n+1}} \boldsymbol{D}\left(\mathbf{u}^{n+1}\right)|^{2}-r|\sqrt{\eta^{n}}\boldsymbol{D}\left(\mathbf{u}^{n}\right)|^{2}\bigg)drdz\\
=&-\left(\int_{\partial\Omega_{L}}2r\eta^{n+1} v_{z}^{n+1}\frac{\partial v_{z}^{n+1}}{\partial z}ds-\int_{\partial\Omega_{L}}2r\eta^{n} v_{z}^{n}\frac{\partial v_{z}^{n}}{\partial z}ds\right)-\frac{1}{4}\int_{\Omega}\bigg(r |\sqrt{\eta^{n+1}} \boldsymbol{D}\left(\mathbf{u}^{n+1}\right)|^{2}-r|\sqrt{\eta^{n}}\boldsymbol{D}\left(\mathbf{u}^{n}\right)|^{2}\bigg)drdz\\
\end{aligned}
$$
We multiply \eqref{37} with $\frac{1}{\alpha}\boldsymbol{R}^{n+1} _t$ to get
\begin{equation}
    \begin{split}
&\frac{1}{\alpha}\frac{\left|\boldsymbol{R}^{n+1}\right|^2}{2}-\frac{1}{\alpha}\frac{\left|\boldsymbol{R}^n\right|^2}{2}+\frac{1}{\alpha}\frac{\left|\boldsymbol{R}^{n+1}-\boldsymbol{R}^n\right|^2}{2}=-\delta_{t}\boldsymbol{R}^{n+1}\int_{\Omega}\nabla\cdot\left(r\eta^{n+1}D(\mathbf{u}^{n})\right)\cdot\mathbf{u}^{n+1}drdz\\
&+\delta_{t}\boldsymbol{R}^{n+1}\int_{\Omega}\left(\mathcal{R}er\rho^{n+1}(\mathbf{u}^{n}\cdot\nabla) \mathbf{u}^{n} \cdot \mathbf{u}^{n+1}+\frac{\mathcal{R}e}{2}\nabla\cdot(r\rho^{n+1}\mathbf{u}^{n})\mathbf{u}^{n} \cdot \mathbf{u}^{n+1}\right) drdz\\
&+\delta_{t}\boldsymbol{R}^{n+1}\int_{\Omega}r\nabla(2p^{n}-p^{n-1})\cdot \mathbf{u}^{n+1}drdz-\delta_{t}\frac{|\boldsymbol{R}^{n+1}|^{2}}{2}\int_{\Omega}r\eta^{n+1}|D(\mathbf{u}^{n})|^{2}drdz\\
&\boldsymbol{R}^{n+1}_{t}\int_{\Omega}\frac{1}{2}\nabla\cdot r\mathbf{J}(\mu^{n})\mathbf{u^{n}}\cdot\mathbf{u}^{n+1}+r\mathbf{J}(\mu^{n})\cdot\nabla\mathbf{u}^{n}\cdot\mathbf{u}^{n+1}drdz+\frac{\delta_{t}\boldsymbol{R}^{n+1}\boldsymbol{K}^{n+1}\mathcal{K}^{n}}{\sqrt{\int_{0}^{t_{n}}-\mathcal{K}d\tau+\boldsymbol{G}}}.\\
    \end{split}
\end{equation}
We multiply \eqref{38} with $2 \delta_t \boldsymbol{K}^{n+1}$ to get
\begin{align}
& \left|\boldsymbol{K}^{n+1}\right|^2-\left|\boldsymbol{K}^n\right|^2+\left|\boldsymbol{K}^{n+1}-\boldsymbol{K}^n\right|^2=\frac{-\delta_t \boldsymbol{K}^{n+1}\boldsymbol{R}^{n+1}\mathcal{K}^{n}}{\sqrt{\int_{0}^{t_{n}}-\mathcal{K}d\tau+\boldsymbol{G}}}.
\end{align}
We multiply \eqref{39} with $\frac{\delta_{t}}{\alpha}\boldsymbol{T}^{n+1}$ to get
\begin{equation}
\frac{1}{\alpha}\frac{\left|\boldsymbol{T}^{n+1}\right|^2-\left|\boldsymbol{T}^n\right|^2+\left|\boldsymbol{T}^{n+1}-\boldsymbol{T}^n\right|^2}{2}=\delta_{t}\boldsymbol{T}^{n+1} \int_{\Omega}\left(\nabla\cdot(r\mathbf{u}^{n+1})\right) p^{n+1} drdz.
\end{equation}
We multiply \eqref{310} with $\frac{\delta_{t}^{2}}{\chi\mathcal{R}e}p^{n+1}$ to get
\begin{align}
&\frac{\delta_{t}^{2}}{2\chi\mathcal{R}e} \int_{\Omega}\big( r\left|\nabla p^{n+1}\right|^2-r\left|\nabla p^n\right|^2+r\left|\nabla p^{n+1}-\nabla p^n\right|^2 \big)drdz=-\delta_{t}\boldsymbol{T}^{n+1} \int_{\Omega}\left(\nabla\cdot(r\mathbf{u}^{n+1})\right) p^{n+1}\big) drdz.
\end{align}
By combining the above equations, we derive
\begin{equation}
    \begin{split}
& \frac{\mathcal{R}e}{2} \int_{\Omega} r \rho^{n+1}\left|\mathbf{u}^{n+1}\right|^2 drdz-\frac{\mathcal{R}e}{2} \int_{\Omega} r \rho^n\left|\mathbf{u}^n\right|^2 drdz+\frac{\mathcal{B}}{2} \epsilon \int_{\Omega}\left(r\left|\nabla \phi^{n+1}\right|^2-r\left|\nabla \phi^n\right|^2\right)drdz \\
& +\frac{\mathcal{B}\boldsymbol{s}}{2 \epsilon} \int_{\Omega}\left(r\left|\phi^{n+1}\right|^2-r\left|\phi^n\right|^2\right) drdz+2\delta_{t}\int_{\partial\Omega_{L}}r\eta^{n+1} v^{n+1}_{z}\frac{\partial v^{n+1}_{z}}{\partial z}ds-2\delta_{t}\int_{\partial\Omega_{L}}r\eta^{n} v^{n}_{z}\frac{\partial v^{n}_{z}}{\partial z}ds \\
&+\frac{\delta_{t}}{4}\int_{\Omega}r\eta^{n+1}|D(\mathbf{u}^{n+1})|^{2}drdz-\frac{\delta_{t}}{4}\int_{\Omega}r\eta^{n}|D(\mathbf{u}^{n})|^{2}drdz+\frac{\delta_{t}^{2}}{2\chi\mathcal{R}e} \int_{\Omega}\big( r\left|\nabla p^{n+1}\right|^2-r\left|\nabla p^n\right|^2\big)drdz\\
&+\frac{1}{\alpha}\frac{\mathcal{B}}{2}\left(\left|\boldsymbol{Q}^{n+1}\right|^2-\left|\boldsymbol{Q}^n\right|^2\right)+\frac{\mathcal{B}}{\epsilon}\left(\left|\boldsymbol{U}^{n+1}\right|^2-\left|\boldsymbol{U}^n\right|^2\right)+\frac{1}{\alpha}\big(\frac{\left|\boldsymbol{R}^{n+1}\right|^2}{2}-\frac{\left|\boldsymbol{R}^n\right|^2}{2}\big)\\
&+\frac{1}{\alpha}\big(\frac{\left|\boldsymbol{T}^{n+1}\right|^2}{2}-\frac{\left|\boldsymbol{T}^n\right|^2}{2}\big)+\left|\boldsymbol{K}^{n+1}\right|^2-\left|\boldsymbol{K}^n\right|^2\leq-\delta_t \mathcal{B} \mathcal{L}_d \mathcal{R}e \int_{\Omega}\left|\nabla r \mu^{n+1}\right|^2 drdz\\
&-\frac{|\boldsymbol{R^{n+1}}|^{2}_t}{2} \int_{\Omega} r \eta^{n+1}\left|D\left(\mathbf{u}^{n}\right)\right|^2 drdz-\delta_{t}\int_{\Omega} 2\eta^{n+1}\frac{(v_{r}^{n+1})^{2}}{r}drdz.
    \end{split}
\end{equation}
Then,
\begin{equation}
\mathbf{E}^{n+1}_{M}-\mathbf{E}^{n}_{M}<-\delta_t \mathcal{B} \mathcal{L}_d\int_{\Omega}\left|\nabla r \mu^{n+1}\right|^2 drdz-\frac{|\boldsymbol{R^{n+1}}|^{2}_t}{2} \int_{\Omega} r \eta^{n+1}\left|D\left(\mathbf{u}^{n}\right)\right|^2 drdz-\delta_{t}\int_{\Omega} 2\eta^{n+1}\frac{(v_{r}^{n+1})^{2}}{r}drdz.
\end{equation}

\subsection{Implementation process and solvability}

In this subsection, we will discuss how to implement $\boldsymbol{STEP 1}$, $\boldsymbol{STEP 2}$ and $\boldsymbol{STEP 3}$. Although we divided the scheme into three steps, it is not the fully decoupling format as expected. Instead, it looks as a coupled scheme since all unknowns are coupled together (e.g. $\boldsymbol{STEP 1}$ couples $\phi^{n+1}$, $\tilde{\mathbf{u}}^{n+1}$,  $\boldsymbol{Q}^{n+1}$, $\boldsymbol{U}^{n+1}$, and $\boldsymbol{STEP 2}$ couples $\mathbf{u}^{n+1}$, $\boldsymbol{R}^{n+1}$, and $\boldsymbol{STEP 3}$ couples $\mathbf{p}^{n+1}$, $\boldsymbol{T}^{n+1}$). Moreover, $\boldsymbol{STEP 1}$, $\boldsymbol{STEP 2}$ and $\boldsymbol{STEP 3}$ also involve many nonlocal terms, which may result in high computational costs. Therefore, to perform  computations, we exploit the nonlocal property of the auxiliary variables  $\boldsymbol{Q}$, $\boldsymbol{U}$, $\boldsymbol{R}$ and  $\boldsymbol{T}$ to get the decoupling implementations and eliminate all nonlocal terms through the following steps.

Based on $\boldsymbol{Q}^{n+1}$, we split the we split $\phi^{n+1}$, $\mu^{n+1}$, $\tilde{\mathbf{u}}^{n+1}$, $U^{n+1}$ to the following linear combination form
\begin{align}
&\left\{\begin{array}{l}
\phi^{n+1}=\phi_{1}^{n+1}+\boldsymbol{Q}^{n+1} \phi_{2}^{n+1},\quad \mu^{n+1}=\mu_{1}^{n+1}+\boldsymbol{Q}^{n+1} \mu_{2}^{n+1}, \label{53} \\ 
\tilde{\mathbf{u}}^{n+1}=\tilde{\mathbf{u}}_{1}^{n+1}+\boldsymbol{Q}^{n+1}\tilde{\mathbf{u}}_{2}^{n+1}, \quad \boldsymbol{U}^{n+1}=\boldsymbol{U}_{1}^{n+1}+\boldsymbol{Q}^{n+1} \boldsymbol{U}_{2}^{n+1}.
\end{array}\right.
\end{align}
The system \eqref{31}, \eqref{32} and \eqref{33} is splitted into two sub-systems as
\begin{align}
&\left\{\begin{array}{l}
r\frac{\phi_{1}^{n+1}}{ t}-\mathcal{L}_{d}\nabla\cdot(r\nabla \mu_{1}^{n+1})=r\frac{\phi^{n}}{ t}, \label{311}\\
r\mu_{1}^{n+1}=-\epsilon \nabla\cdot(r\nabla\phi_{1}^{n+1})+\frac{r}{\epsilon} H^{n} \boldsymbol{U}_{1}^{n+1}+\frac{rs}{\epsilon} \phi_{1}^{n+1}, \\
\tilde{\mathbf{u}}_{1}^{n+1}=\mathbf{u}^{n},
\end{array}\right. \\
&\left\{\begin{array}{l}
r\frac{\phi_{2}^{n+1}}{ t}-\mathcal{L}_{d}\nabla\cdot(r\nabla \mu_{2}^{n+1})=-r\mathbf{u}^{n} \cdot\nabla \phi^{n}, \\
r\mu_{2}^{n+1}=-\epsilon \nabla\cdot(r\nabla\phi_{2}^{n+1})+\frac{r}{\epsilon} H^{n} \boldsymbol{U}_{2}^{n+1}+\frac{rs}{\epsilon} \phi_{2}^{n+1}, \label{312}\\
\rho^{n} \frac{\tilde{\mathbf{u}}_{2}^{n+1}}{ t}- \mathcal{B}\mu^{n}\nabla\phi^{n}=0.
\end{array}\right.
\end{align}
To solve the two subsystems \eqref{311} and \eqref{312}, we continue to employ the splitting technique, where the variables $\phi_{1}^{n+1}$, $\phi_{2}^{n+1}$, $\mu_{1}^{n+1}$ and $\mu_{2}^{n+1}$ are decomposed into a linear combination form involving the nonlocal variables $\boldsymbol{U}^{n+1}_{1}$, $\boldsymbol{U}^{n+1}_{2}$, respectively, which reads as
\begin{align}
&\left\{\begin{array}{l}
\phi_{1}^{n+1}=\phi_{11}^{n+1}+\boldsymbol{U}_{1}^{n+1} \phi_{12}^{n+1}, \quad \mu_{1}^{n+1}=\mu_{11}^{n+1}+\boldsymbol{U}_{1}^{n+1} \mu_{12}^{n+1}, \label{56}\\ 
\phi_{2}^{n+1}=\phi_{21}^{n+1}+\boldsymbol{U}_{2}^{n+1} \phi_{22}^{n+1}, \quad \mu_{2}^{n+1}=\mu_{21}^{n+1}+\boldsymbol{U}_{2}^{n+1} \mu_{22}^{n+1}.
\end{array}\right.
\end{align}
By utilizing \eqref{56} to substitute $\phi_{1}^{n+1}$, $\phi_{2}^{n+1}$, $\mu_{1}^{n+1}$, $\mu_{2}^{n+1}$ in \eqref{311} and \eqref{312}, and then splitting the results based on $\boldsymbol{U}^{n+1}_{1}$ and $\boldsymbol{U}^{n+1}_{2}$, respectively, we get
\begin{align}
&\left\{\begin{array}{l}
r\frac{\phi_{11}^{n+1}}{ t}-\mathcal{L}_{d}\nabla\cdot(r\nabla \mu_{11}^{n+1})=r\frac{\phi^{n}}{ t}, \label{313} \\
r\mu_{11}^{n+1}=-\epsilon \nabla\cdot(r\nabla\phi_{11}^{n+1})+\frac{rs}{\epsilon} \phi_{11}^{n+1},
\end{array}\right. \\
&\left\{\begin{array}{l}
r\frac{\phi_{12}^{n+1}}{ t}-\mathcal{L}_{d} \nabla\cdot(r\nabla \mu_{12}^{n+1})=0,\label{314}\\
r\mu_{12}^{n+1}=-\epsilon  \nabla\cdot(r\nabla\phi_{12}^{n+1})+\frac{rs}{\epsilon} \phi_{12}^{n+1}+\frac{r}{\epsilon} H^{n}, 
\end{array}\right. \\
&\left\{\begin{array}{l}
r\frac{\phi_{21}^{n+1}}{ t}-\mathcal{L}_{d} \nabla\cdot(r\nabla \mu_{21}^{n+1})=-r\mathbf{u}^{n} \cdot\nabla \phi^{n}, \label{315}\\
r\mu_{21}^{n+1}=-\epsilon  \nabla\cdot(r\nabla\phi_{21}^{n+1})+\frac{rs}{\epsilon} \phi_{21}^{n+1},
\end{array}\right. \\
&\left\{\begin{array}{l}
r\frac{\phi_{22}^{n+1}}{ t}-\mathcal{L}_{d} \nabla\cdot(r\nabla \mu_{22}^{n+1})=0, \label{316}\\
r\mu_{22}^{n+1}=-\epsilon  \nabla\cdot(r\nabla\phi_{22}^{n+1})+\frac{rs}{\epsilon} \phi_{22}^{n+1}+\frac{r}{\epsilon} H^{n}.
\end{array}\right.
\end{align}
The boundary conditions for \eqref{313}-\eqref{316} are 
\begin{align}
& \phi_{11}^{n+1}|_{\partial\Omega_{L}}=\phi^{n+1}|_{\partial\Omega_{L}}, \phi_{12}^{n+1}|_{\partial\Omega_{L}}=\phi_{21}^{n+1}|_{\partial\Omega_{L}}=\phi_{22}^{n+1}|_{\partial\Omega_{L}}=0,\\
& \partial_{\mathbf{n}}\phi_{11}^{n+1}|_{\partial\Omega\setminus\partial\Omega_{L}}=\partial_{\mathbf{n}}\phi_{12}^{n+1}|_{\partial\Omega\setminus\partial\Omega_{L}}=\partial_{\mathbf{n}}\phi_{21}^{n+1}|_{\partial\Omega\setminus\partial\Omega_{L}}=\partial_{\mathbf{n}}\phi_{22}^{n+1}|_{\partial\Omega\setminus\partial\Omega_{L}}=0,\\
& \mu_{11}^{n+1}|_{\partial\Omega_{L}}=\mu_{12}^{n+1}|_{\partial\Omega_{L}}=\mu_{21}^{n+1}|_{\partial\Omega_{L}}=\mu_{22}^{n+1}|_{\partial\Omega_{L}}=0,\\
& \partial_{\mathbf{n}}\mu_{11}^{n+1}|_{\partial\Omega\setminus\partial\Omega_{L}}=\partial_{\mathbf{n}}\mu_{12}^{n+1}|_{\partial\Omega\setminus\partial\Omega_{L}}=\partial_{\mathbf{n}}\mu_{21}^{n+1}|_{\partial\Omega\setminus\partial\Omega_{L}}=\partial_{\mathbf{n}}\mu_{22}^{n+1}|_{\partial\Omega\setminus\partial\Omega_{L}}=0.
\end{align}
Then we will solve $\boldsymbol{U}^{n+1}_{1}$ and $\boldsymbol{U}^{n+1}_{2}$. By employing the split form of $\boldsymbol{U}^{n+1}$ and $\phi^{n+1}$ in \eqref{53}, we can obtain
\begin{align}
\left\{\begin{array}{l}
\boldsymbol{U}_{1}^{n+1}=\frac{1}{2} \int_{\Omega} rH^{n} \phi_{1}^{n+1} drdz+g^{n}, \label{317}\\
\boldsymbol{U}_{2}^{n+1}=\frac{1}{2} \int_{\Omega} rH^{n} \phi_{2}^{n+1} drdz.
\end{array}\right.
\end{align}
where, $g^{n}=\boldsymbol{U}^{n}-\frac{1}{2} \int_{\Omega} rH^{n} \phi^{n} drdz$. After applying a simple factorization given in \eqref{56}, we can get
\begin{align}
\left\{\begin{array}{l}
\boldsymbol{U}_{1}^{n+1}=\frac{\frac{1}{2} \int_{\Omega}r H^{n} \phi_{11}^{n+1} drdz+g^{n}}{1-\frac{1}{2} \int_{\Omega}r H^{n} \phi_{12}^{n+1} drdz}, \\
\boldsymbol{U}_{2}^{n+1}=\frac{\frac{1}{2} \int_{\Omega}r H^{n} \phi_{21}^{n+1} drdz}{1-\frac{1}{2} \int_{\Omega}r H^{n} \phi_{22}^{n+1} drdz}.
\end{array}\right.
\end{align}
We demonstrate that $\boldsymbol{U}^{n+1}_{1}$ and $\boldsymbol{U}^{n+1}_{2}$ are solvable by verifying the denominators are non-zero. This can be achieved by applying a simple energy estimate to the subsystem \eqref{314}. For the first equation in \eqref{314}, we take the inner product with $\delta_{t}\mu^{n+1}_{12}$ in the $L^{2}$ space, then we have
\begin{equation}
\int_{\Omega}r\phi_{12}^{n+1}\mu^{n+1}_{12}drdz+\delta_{t}\mathcal{L}_{d}\int_{\Omega}r|\nabla \mu_{12}^{n+1}|^{2}drdz=0.\label{U_1}
\end{equation}
For the second equation in \eqref{314}, we take the inner product with $\phi_{12}^{n+1}$ in the $L^{2}$ space,
\begin{equation}
\int_{\Omega}r\mu_{12}^{n+1}\phi^{n+1}_{12}drdz=\epsilon \int_{\Omega} r|\nabla\phi_{12}^{n+1}|^{2}drdz+\frac{s}{\epsilon}\int_{\Omega}r|\phi_{12}^{n+1}|^{2}drdz+\frac{1}{\epsilon}\int_{\Omega}rH^{n}\phi^{n+1}_{12}drdz.\label{U_2}
\end{equation}
By combining the above two equations \eqref{U_1} and \eqref{U_2}, we have
\begin{equation}
\frac{1}{2}\int_{\Omega}rH^{n}\phi^{n+1}_{12}drdz=-\delta_{t}\frac{\epsilon}{2}\mathcal{L}_{d}\int_{\Omega}r|\nabla \mu_{12}^{n+1}|^{2}drdz-\frac{\epsilon^{2}}{2}\int_{\Omega} r|\nabla\phi_{12}^{n+1}|^{2}drdz-\frac{s}{2}\int_{\Omega}r|\phi_{12}^{n+1}|^{2}drdz.
\end{equation}
Then, we get
\begin{equation}
1-\frac{1}{2} \int_{\Omega}r H^{n} \phi_{12}^{n+1} drdz=1+\delta_{t}\frac{\epsilon}{2}\mathcal{L}_{d}\int_{\Omega}r|\nabla \mu_{12}^{n+1}|^{2}drdz+\frac{\epsilon^{2}}{2}\int_{\Omega} r|\nabla\phi_{12}^{n+1}|^{2}drdz+\frac{s}{2}\int_{\Omega}r|\phi_{12}^{n+1}|^{2}drdz> 0.
\end{equation}
Next, we solve $\boldsymbol{Q}^{n+1}$ from \eqref{35}. We use the split form of $\mu^{n+1}$, $\tilde{\mathbf{u}}^{n+1}$ to rewrite \eqref{35} to the following equation,
\begin{equation}
\frac{\boldsymbol{Q}^{n+1}-\boldsymbol{Q}^{n}}{\delta_{t}}=\alpha\int_{\Omega}\big(r\mathbf{u}^{n}\cdot\nabla \phi^{n} (\mu_{1}^{n+1}+\boldsymbol{Q}^{n+1} \mu_{2}^{n+1})-r(\tilde{\mathbf{u}}_{1}^{n+1}+\boldsymbol{Q}^{n+1}\tilde{\mathbf{u}}_{2}^{n+1})\cdot\nabla \phi^{n} \mu^{n}\big) drdz.
\end{equation}
Then, we obtain
\begin{equation}
\boldsymbol{Q}^{n+1}\left(\frac{1}{\delta_{t}}-\alpha\int_{\Omega}\left(r\mathbf{u}^{n}\cdot\nabla\phi^{n}\mu_{2}^{n+1}-r\tilde{\mathbf{u}}_{2}^{n+1}\cdot\nabla\phi^{n}\mu^{n}\right)drdz\right)=\frac{\boldsymbol{Q}^{n}}{\delta_{t}}+\alpha\int_{\Omega}\left(r\mathbf{u}^{n}\cdot\nabla\phi^{n}\mu_{1}^{n+1}-r\tilde{\mathbf{u}}_{1}^{n+1}\cdot\nabla\phi^{n}\mu^{n}\right)drdz.
\end{equation}
We need to verify that $\boldsymbol{Q}^{n+1}$ is solvable by showing 
\begin{equation}
\frac{1}{\delta_{t}}-\alpha\int_{\Omega}\big(r\mathbf{u}^{n}\cdot\nabla\phi^{n}\mu_{2}^{n+1}-r\tilde{\mathbf{u}}_{2}^{n+1}\cdot\nabla\phi^{n}\mu^{n}\big)drdz \neq 0.
\end{equation}
For $\mathcal{R}e\rho^{n} \frac{\tilde{\mathbf{u}}_{2}^{n+1}}{ t}-\mathcal{B}\mu^{n}\nabla\phi^{n}=0$ in \eqref{312}, we take the $L^{2}$ inner product of it with $\nabla\phi^{n}\mu^{n}$, i.e.,
\begin{equation}
\int_{\Omega}\left(r\tilde{\mathbf{u}}_{2}^{n+1}\cdot\nabla\phi^{n}\mu^{n}\right)drdz=\int_{\Omega}\left(r\frac{\mathcal{B}_{t}}{\rho_{n}\mathcal{R}e}|\nabla\phi^{n}\mu^{n}|^{2}\right)drdz.
\end{equation}
By taking the $L^{2}$ inner product of the first equation in \eqref{312} with $\delta_{t}\mu_{2}^{n+1}$, of the second equation of \eqref{312} with $\phi_{2}^{n+1}$, and combining the obtained two equations, we get
\begin{equation}
-\int_{\Omega}r\mathbf{u}^{n}\cdot\nabla\phi^{n}\mu_{2}^{n+1}drdz=\mathcal{L}_{d}\int_{\Omega}r|\nabla\mu^{n+1}|^{2}drdz+\frac{\epsilon}{\delta_{t}}\int_{\Omega}r|\nabla\phi_{2}^{n+1}|^{2}drdz+\frac{1}{\delta_{t}\epsilon}|U^{n+1}_{2}|^{2}+\frac{s}{\epsilon\delta_{t}}\int_{\Omega}r|\phi_{2}^{n+1}|^2drdz,
\end{equation}
which implies
\begin{equation}
\begin{split}
&\frac{1}{\delta_{t}}-\alpha\int_{\Omega}\left(\mathbf{u}^{n}\cdot\nabla\phi^{n}\mu_{2}^{n+1}-\tilde{\mathbf{u}}_{2}^{n+1}\cdot\nabla\phi^{n}\mu^{n}\right)drdz=\frac{1}{\delta_{t}}+\alpha\mathcal{L}_{d}\int_{\Omega}r|\nabla\mu^{n+1}|^{2}drdz+\alpha\frac{\epsilon}{\delta_{t}}\int_{\Omega}r|\nabla\phi_{2}^{n+1}|^{2}drdz\\
&+\frac{\alpha}{\delta_{t}\epsilon}|U^{n+1}_{2}|^{2}+\frac{s\alpha}{\epsilon\delta_{t}}\int_{\Omega}r|\phi_{2}^{n+1}|^2drdz+\alpha\int_{\Omega}\left(r\frac{\delta_{t}}{\rho_{n}\mathcal{R}e}|\nabla\phi^{n}\mu^{n}|^{2}\right)drdz> 0.
\end{split}
\end{equation}
Next, we solve $\boldsymbol{R}^{n+1}$ from \eqref{37}. Using the nonlocal variable $\boldsymbol{R}^{n+1}$, we can rewrite $\mathbf{u}^{n+1}$ and $\boldsymbol{K}^{n+1}$ to be the following linear form as
\begin{equation}
\mathbf{u}^{n+1}=\mathbf{u}^{n+1}_{1}+\boldsymbol{R}^{n+1}\mathbf{u}^{n+1}_{2}, \boldsymbol{K}^{n+1}=\boldsymbol{K}^{n+1}_{1}+\boldsymbol{R}^{n+1}\boldsymbol{K}^{n+1}_{2}.\label{R_1}
\end{equation}
Using \eqref{R_1}, we decompose the equation \eqref{36} into the following two sub-equations
according to $\boldsymbol{R}^{n+1}$,
\begin{equation}
    \begin{split}
\mathcal{R}e\frac{r}{2}\frac{\rho^{n+1}+\rho^{n}}{\delta_{t}}\mathbf{u}_{1}^{n+1}-\nabla\cdot\left(r\eta^{n+1}D(\mathbf{u}_{1}^{n+1})\right)+2\frac{\eta^{n+1}}{r}(0,v^{n+1}_{1_{r}})=&\mathcal{R}er\rho^{n}\frac{\tilde{\mathbf{u}}^{n+1}}{\delta_{t}} -\nabla\cdot\left(r\sqrt{\eta^{n+1}\eta^{n}}D(\mathbf{u}^{n})\right),
 \end{split}
\end{equation}
\begin{equation}
    \begin{split}
&\mathcal{R}e\frac{r}{2}\frac{\rho^{n+1}+\rho^{n}}{\delta_{t}}\mathbf{u}_{2}^{n+1}-\nabla\cdot\left(r\eta^{n+1}D(\mathbf{u}_{2}^{n+1})\right)+2\frac{\eta^{n+1}}{r}(0,v^{n+1}_{2_{r}})=-\mathcal{R}e\rho^{n}\left(r\mathbf{u}^{n}\cdot\nabla\right) \mathbf{u}^{n} \\ 
&-\frac{\mathcal{R}e}{2}\left(\nabla\cdot(r\rho^{n+1}\mathbf{u}^{n})\right)\mathbf{u}^{n}+\nabla\cdot\left(r\eta^{n+1}D(\mathbf{u}^{n})\right)-r\nabla(2p^{n}-p^{n-1})-\frac{1}{2}\nabla\cdot r\mathbf{J}(\mu^{n})\mathbf{u}^{n}-r\mathbf{J}(\mu^{n})\cdot\nabla\mathbf{u}^{n}.    
\end{split}\label{63}
\end{equation}
The boundary conditions of the above two equations are 
\begin{align}
&{\mathbf{u}_{1}^{n+1}}|_{\Gamma_2\cup\Gamma_3}=\mathbf{u}^{n+1}|_{\Gamma_2\cup\Gamma_3}, {\mathbf{u}_{2}}^{n+1}|_{\Gamma_2\cup\Gamma_3}=0,\mathbf{u}_{1}^{n+1}|_{\Gamma_1}=\mathbf{u}_{2}^{n+1}|_{\Gamma_1}=0,\\
&{v^{n+1}_{1_{r}}}|_{\Gamma_4\cup\Gamma_5}={v^{n+1}_{2_{r}}}|_{\Gamma_4\cup\Gamma_5}=0, \partial_{\mathbf{n}}{v^{n+1}_{1_{z}}}|_{\Gamma_4\cup\Gamma_5}=\partial_{\mathbf{n}}{v^{n+1}_{2_{z}}}|_{\Gamma_4\cup\Gamma_5}=0.
\end{align}
Substituting the linear combination form of $\boldsymbol{K}^{n+1}$, we obtain
\begin{equation}
    \begin{split}
\frac{\boldsymbol{R}^{n+1}-\boldsymbol{R}^{n}}{\delta_{t}}=\alpha\bigg(&\int_{\Omega}\big(\mathcal{R}er\rho^{n+1}(\mathbf{u}^{n}\cdot\nabla) \mathbf{u}^{n} \cdot (\mathbf{u}^{n+1}_{1}+\boldsymbol{R}^{n+1}\mathbf{u}^{n+1}_{2})+\frac{\mathcal{R}e}{2}\nabla\cdot(r\rho^{n+1}\mathbf{u}^{n})\mathbf{u}^{n} \cdot (\mathbf{u}^{n+1}_{1}+\boldsymbol{R}^{n+1}\mathbf{u}^{n+1}_{2})\big) drdz\\
&+\int_{\Omega}r\nabla(2p^{n}-p^{n-1})\cdot(\mathbf{u}^{n+1}_{1}+\boldsymbol{R}^{n+1}\mathbf{u}^{n+1}_{2})drdz-\int_{\Omega}\nabla\cdot\left(r\eta^{n+1}D(\mathbf{u}^{n})\right)\cdot(\mathbf{u}^{n+1}_{1}+\boldsymbol{R}^{n+1}\mathbf{u}^{n+1}_{2})drdz\\
&+\int_{\Omega}\frac{1}{2}\nabla\cdot r\mathbf{J}(\mu^{n})\mathbf{u}^{n}\cdot(\mathbf{u}^{n+1}_{1}+\boldsymbol{R}^{n+1}\mathbf{u}^{n+1}_{2})+r\mathbf{J}(\mu^{n})\cdot\nabla\mathbf{u}^{n}\cdot(\mathbf{u}^{n+1}_{1}+\boldsymbol{R}^{n+1}\mathbf{u}^{n+1}_{2})drdz\\
&-\frac{\boldsymbol{R}^{n+1}}{2}\int_{\Omega}r\eta^{n+1}|D(\mathbf{u}^{n})|^{2}drdz+\frac{(\boldsymbol{K}^{n+1}_{1}+\boldsymbol{R}^{n+1}\boldsymbol{K}^{n+1}_{2})\mathcal{K}^{n}}{\sqrt{-\int_{0}^{t_{n}}\mathcal{K}d\tau+\boldsymbol{G}}}\bigg).\\
    \end{split}\label{R_2}
\end{equation}
Then, we rewrite \eqref{R_2} as,
\begin{equation}
    \begin{split}
&\boldsymbol{R}^{n+1}\bigg(\frac{1}{\delta_{t}}-\alpha\int_{\Omega}\left(\mathcal{R}er\rho^{n+1}(\mathbf{u}^{n}\cdot\nabla) \mathbf{u}^{n}\cdot\mathbf{u}^{n+1}_{2}+\frac{\mathcal{R}e}{2}\nabla\cdot(r\rho^{n+1}\mathbf{u}^{n})\mathbf{u}^{n} \cdot\mathbf{u}^{n+1}_{2}\right)drdz  -\alpha\int_{\Omega}r\nabla(2p^{n}-p^{n-1})\cdot\mathbf{u}^{n+1}_{2}drdz\\ 
&+\alpha\int_{\Omega}\nabla\cdot\left(r\eta^{n+1}D(\mathbf{u}^{n})\right)\cdot\mathbf{u}^{n+1}_{2}drdz+\frac{\alpha}{2}\int_{\Omega}r\eta^{n+1}|D(\mathbf{u}^{n})|^{2}drdz-\alpha\frac{\boldsymbol{K}^{n+1}_{2}\mathcal{K}^{n}}{\sqrt{-\int_{0}^{t_{n}}\mathcal{K}d\tau+\boldsymbol{G}}}-\alpha\int_{\Omega}\frac{1}{2}\nabla\cdot r\mathbf{J}(\mu^{n})\mathbf{u}^{n}\cdot\mathbf{u}^{n+1}_{2}\\
&+r\mathbf{J}(\mu^{n})\cdot\nabla\mathbf{u}^{n}\cdot\mathbf{u}^{n+1}_{2}drdz\bigg)\\
=&\frac{\boldsymbol{R}^{n}}{\delta_{t}}+\alpha\int_{\Omega}\left(\mathcal{R}er\rho^{n+1}(\mathbf{u}^{n}\cdot\nabla) \mathbf{u}^{n} \cdot \mathbf{u}^{n+1}_{1}+\frac{\mathcal{R}e}{2}\nabla\cdot(r\rho^{n+1}\mathbf{u}^{n})\mathbf{u}^{n} \cdot \mathbf{u}^{n+1}_{1}\right) drdz-\alpha\int_{\Omega}\nabla\cdot\left(r\eta^{n+1}D(\mathbf{u}^{n})\right)\cdot\mathbf{u}^{n+1}_{1}drdz\\
&+\alpha\int_{\Omega}r\nabla(2p^{n}-p^{n-1})\cdot \mathbf{u}^{n+1}_{1}drdz+\alpha\frac{\boldsymbol{K}^{n+1}_{1}\mathcal{K}^{n}}{\sqrt{-\int_{0}^{t_{n}}\mathcal{K}d\tau+\boldsymbol{G}}}+\alpha\int_{\Omega}\frac{1}{2}\nabla\cdot r\mathbf{J}(\mu^{n})\mathbf{u}^{n}\cdot\mathbf{u}^{n+1}_{1}+r\mathbf{J}(\mu^{n})\cdot\nabla\mathbf{u}^{n}\cdot\mathbf{u}^{n+1}_{1}drdz.
    \end{split}
\end{equation}
We need to verify that $R^{n+1}$ is solvable by showing that 
\begin{equation}
    \begin{split}
&\frac{1}{\delta_{t}}-\alpha\int_{\Omega}\left(\mathcal{R}er\rho^{n+1}(\mathbf{u}^{n}\cdot\nabla) \mathbf{u}^{n}\cdot\mathbf{u}^{n+1}_{2}+\frac{\mathcal{R}e}{2}\nabla\cdot(r\rho^{n+1}\mathbf{u}^{n})\mathbf{u}^{n} \cdot\mathbf{u}^{n+1}_{2}\right)drdz   -\alpha\int_{\Omega}r\nabla(2p^{n}-p^{n-1})\cdot\mathbf{u}^{n+1}_{2}drdz \\ 
&+\alpha\int_{\Omega}\nabla\cdot\left(r\eta^{n+1}D(\mathbf{u}^{n})\right)\cdot\mathbf{u}^{n+1}_{2}drdz+\frac{\alpha}{2}\int_{\Omega}r\eta^{n+1}|D(\mathbf{u}^{n})|^{2}drdz-\alpha\frac{\boldsymbol{K}^{n+1}_{2}\mathcal{K}^{n}}{\sqrt{-\int_{0}^{t_{n}}\mathcal{K}d\tau+\boldsymbol{G}}}-\alpha\int_{\Omega}\frac{1}{2}\nabla\cdot r\mathbf{J}(\mu^{n})\mathbf{u}^{n}\cdot\mathbf{u}^{n+1}_{2}\\
&+r\mathbf{J}(\mu^{n})\cdot\nabla\mathbf{u}^{n}\cdot\mathbf{u}^{n+1}_{2}drdz\bigg)\neq 0
    \end{split}\label{63_3}
\end{equation}
We take the inner product of \eqref{63} with $\mathbf{u}^{n+1}_{2}$, then we have 
\begin{equation}
    \begin{split}
&-\int_{\Omega}\bigg(\mathcal{R}e\rho^{n}(r\mathbf{u}^{n}\cdot\nabla)\mathbf{u}^{n}\cdot\mathbf{u}_{2}^{n+1}+\frac{\mathcal{R}e}{2}\nabla\cdot(r\rho^{n}\mathbf{u}^{n})\mathbf{u}^{n}\cdot\mathbf{u}_{2}^{n+1}+r\nabla(2p^{n}-p^{n-1})\cdot\mathbf{u}_{2}^{n+1}\bigg)drdz  \\ 
&+\int_{\Omega}\nabla\cdot\left(r\eta^{n+1}D(\mathbf{u}^{n})\right)\cdot\mathbf{u}^{n+1}_{2}drdz-\int_{\Omega}\frac{1}{2}\nabla\cdot r\mathbf{J}(\mu^{n})\mathbf{u}^{n}\cdot\mathbf{u}^{n+1}_{2}+r\mathbf{J}(\mu^{n})\cdot\nabla\mathbf{u}^{n}\cdot\mathbf{u}^{n+1}_{2}drdz\\ 
=&\int_{\Omega}\frac{\mathcal{R}e}{2  t}\left(r\rho^{n+1}+r\rho^{n}\right) |\mathbf{u}_{2}^{n+1}|^{2}+\frac{1}{2}\int_{\Omega}r\eta^{n+1} |D\left(\mathbf{u}_{2}^{n+1}\right)|^{2}drdz+2\int_{\Omega}\frac{\eta^{n+1}}{r}(v^{n+1}_{2_{r}})^{2}drdz>0. 
    \end{split}\label{63_2}
\end{equation}
Using the split form of the variable $\boldsymbol{K}^{n+1}$, we rewrite \eqref{38} as the following form
\begin{equation}
\frac{\boldsymbol{K}_{1}^{n+1}+\boldsymbol{R}^{n+1}\boldsymbol{K}_{2}^{n+1}-\boldsymbol{K}^{n}}{\delta_{t}}=\frac{\boldsymbol{R}^{n+1}}{2}\frac{-\mathcal{K}^{n}}{\sqrt{-\int_{0}^{t_{n}}\mathcal{K}d\tau+\boldsymbol{G}}}.
\end{equation}
Then, we can get
\begin{align}
&\frac{\boldsymbol{K}_{2}^{n+1}}{\delta_{t}}=\frac{-\mathcal{K}^{n}}{2\sqrt{-\int_{0}^{t_{n}}\mathcal{K}d\tau+\boldsymbol{G}}},\label{64}\\
&\frac{\boldsymbol{K}_{1}^{n+1}-\boldsymbol{K}^{n}}{\delta_{t}}=0.\label{65}
\end{align}
We multiply \eqref{64} with $\boldsymbol{K}_{2}^{n+1}$ to get
\begin{align}
\boldsymbol{K}_{2}^{n+1}\frac{-\mathcal{K}^{n}}{\sqrt{-\int_{0}^{t_{n}}\mathcal{K}d\tau+\boldsymbol{G}}}=\frac{2}{\delta_{t}}|\boldsymbol{K}_{2}^{n+1}|^{2}. \label{66}
\end{align}
By combining \eqref{63_3}, \eqref{63_2} and \eqref{66}, we can obtain
\begin{equation}
    \begin{split}
&\frac{1}{\delta_{t}}-\alpha\int_{\Omega}\left(\mathcal{R}er\rho^{n+1}(\mathbf{u}^{n}\cdot\nabla) \mathbf{u}^{n}\cdot\mathbf{u}^{n+1}_{2}+\frac{\mathcal{R}e}{2}\nabla\cdot(r\rho^{n+1}\mathbf{u}^{n})\mathbf{u}^{n} \cdot\mathbf{u}^{n+1}_{2}\right)drdz  -\alpha\int_{\Omega}r\nabla(2p^{n}-p^{n-1})\cdot\mathbf{u}^{n+1}_{2}drdz\\ 
&+\alpha\int_{\Omega}\nabla\cdot\left(r\eta^{n+1}D(\mathbf{u}^{n})\right)\cdot\mathbf{u}^{n+1}_{2}drdz+2\alpha\int_{\Omega}\frac{\eta^{n+1}}{r}(v^{n+1}_{2_{r}})^{2}drdz-\alpha\frac{\boldsymbol{K}^{n+1}_{2}\mathcal{K}^{n}}{\sqrt{-\int_{0}^{t_{n}}\mathcal{K}d\tau+\boldsymbol{G}}}\\
=&\frac{1}{\delta_{t}}+\frac{2\alpha}{\delta_{t}}|\boldsymbol{K}_{2}^{n+1}|^{2}+\alpha\int_{\Omega}\frac{\mathcal{R}e}{2  t}\left(r\rho^{n+1}+r\rho^{n}\right) |\mathbf{u}_{2}^{n+1}|^{2}drdz+\frac{\alpha}{2}\int_{\Omega}r\eta^{n+1} |D\left(\mathbf{u}_{2}^{n+1}|^{2}\right)drdz+2\alpha\int_{\Omega}\frac{\eta^{n+1}}{r}(v^{n+1}_{2_{r}})^{2}drdz\\
&+\alpha\frac{1}{2}\int_{\Omega}r\eta^{n+1}|D(\mathbf{u}^{n})|^{2}drdz>0.
    \end{split}
\end{equation}
Finally, we will solve $\boldsymbol{T}^{n+1}$ from \eqref{39}, we rewrite $p^{n+1}$ to be a linear form as
\begin{align}
p^{n+1}=p_{1}^{n+1}+\boldsymbol{T}^{n+1}p_{2}^{n+1}. \label{69}
\end{align}
By \eqref{69}, we rewrite \eqref{39} as the following form
\begin{align}
\boldsymbol{T}^{n+1}\left(\frac{1}{\delta_{t}}-\alpha\int_{\Omega}\nabla\cdot(r\mathbf{u}^{n+1})p_{2}^{n+1}drdz\right)=\alpha\int_{\Omega}\nabla\cdot(r\mathbf{u}^{n+1})p_{1}^{n+1}drdz. \label{610}
\end{align}
We need to verify that the following equation
\begin{align}
\frac{1}{\delta_{t}}-\alpha\int_{\Omega}\nabla\cdot(r\mathbf{u}^{n+1})p_{2}^{n+1}drdz\neq 0.
\end{align}
Use \eqref{69} to replace $p^{n+1}$ in \eqref{310}, we obtain
\begin{align}
&p_{1}^{n+1}=p^{n},\label{72}\\
&\nabla\cdot \left(r\nabla p_{2}^{n+1}\right)=\frac{\chi\mathcal{R}e}{\delta_{t}}\nabla\cdot(r\mathbf{u}^{n+1}).\label{73}
\end{align}
The boundary conditions are
\begin{align}
\partial_{\mathbf{n}} p_{2}^{n+1}|_{\partial\Omega\setminus\partial\Omega_{R}}=0,p_{2}^{n+1}|_{\partial\Omega_{R}}=0.
\end{align}
For the equation\eqref{73}, we take the inner product with $p_{2}^{n+1}$ in the $L^{2}$ space, we can get 
\begin{align}
\int_{\Omega}r|\nabla p_{2}^{n+1}|^{2}drdz=-\frac{\chi\mathcal{R}e}{\delta_{t}}\int_{\Omega}\nabla\cdot(r\mathbf{u}^{n+1})p_{2}^{n+1}drdz.\label{74}
\end{align}
By \eqref{74}, we can obtain
\begin{equation}
\frac{1}{\delta_{t}}-\alpha\int_{\Omega}\nabla\cdot(r\mathbf{u}^{n+1})p_{2}^{n+1}drdz=\frac{1}{\delta_{t}}+\frac{\alpha_{t}}{\chi\mathcal{R}e}\int_{\Omega}r|\nabla p_{2}^{n+1}|^{2}drdz>0.
\end{equation}
From the decoupled implementations described above, we can see that the split form of all variables is as follows
\begin{align}
\begin{cases}
\phi^{n+1}=\phi_{1}^{n+1}+\boldsymbol{Q}^{n+1} \phi_{2}^{n+1}, & \mu^{n+1}=\mu_{1}^{n+1}\boldsymbol{Q}^{n+1} \mu_{2}^{n+1}, \\ \tilde{\mathbf{u}}^{n+1}=\tilde{\mathbf{u}}_{1}^{n+1}+\boldsymbol{Q}^{n+1}\tilde{\mathbf{u}}_{2}^{n+1}, & \boldsymbol{U}^{n+1}=\boldsymbol{U}_{1}^{n+1}+\boldsymbol{Q}^{n+1} \boldsymbol{U}_{2}^{n+1},
\end{cases}
\end{align}
\begin{align}
\begin{cases}\phi_{1}^{n+1}=\phi_{11}^{n+1}+\boldsymbol{U}_{1}^{n+1} \phi_{12}^{n+1}, & \mu_{1}^{n+1}=\mu_{11}^{n+1}+\boldsymbol{U}_{1}^{n+1} \mu_{12}^{n+1}, \\ \phi_{2}^{n+1}=\phi_{21}^{n+1}+\boldsymbol{U}_{2}^{n+1} \phi_{22}^{n+1}, & \mu_{2}^{n+1}=\mu_{21}^{n+1}+\boldsymbol{U}_{2}^{n+1} \mu_{22}^{n+1},\end{cases}
\end{align}
\begin{equation}
\mathbf{u}^{n+1}=\mathbf{u}^{n+1}_{1}+\boldsymbol{R}^{n+1}\mathbf{u}^{n+1}_{2}, \boldsymbol{K}^{n+1}=\boldsymbol{K}^{n+1}_{1}+\boldsymbol{R}^{n+1}\boldsymbol{K}^{n+1}_{2},
\end{equation}
\begin{equation}
p^{n+1}=p_{1}^{n+1}+\boldsymbol{T}^{n+1}p_{2}^{n+1}.
\end{equation}
By the above linear combination form, we can rewrite \eqref{31}-\eqref{pb} and decompose them into the following subequations.\\
$\boldsymbol{STEP 1:}$
\begin{align}
& \left\{\begin{array}{l}
r\frac{\phi_{11}^{n+1}-\phi^{n}}{\delta_{t}}-\mathcal{L}_{d} \left[\frac{\partial}{\partial r}\left(r\frac{\partial \mu_{11}^{n+1}}{\partial r}\right)+\frac{\partial}{\partial z}\left(r\frac{\partial \mu_{11}^{n+1}}{\partial z}\right)\right]=0, \\
r\mu_{11}^{n+1}=-\epsilon \left[\frac{\partial}{\partial r}\left(r\frac{\partial \phi_{11}^{n+1}}{\partial r}\right)+\frac{\partial}{\partial z}\left(r\frac{\partial \phi_{11}^{n+1}}{\partial z}\right)\right]+\frac{rs}{\epsilon}\phi_{11}^{n+1}, \\
\end{array}\right. \\
& \left\{\begin{array}{l}
r\frac{\phi_{12}^{n+1}}{\delta_{t}}-\mathcal{L}_{d} \left[\frac{\partial}{\partial r}\left(r\frac{\partial \mu_{12}^{n+1}}{\partial r}\right)+\frac{\partial}{\partial z}\left(r\frac{\partial \mu_{12}^{n+1}}{\partial z}\right)\right]=0, \\
r\mu_{12}^{n+1}=-\epsilon \left[\frac{\partial}{\partial r}\left(r\frac{\partial \phi_{12}^{n+1}}{\partial r}\right)+\frac{\partial}{\partial z}\left(r\frac{\partial \phi_{12}^{n+1}}{\partial z}\right)\right]+\frac{rs}{\epsilon}\phi_{12}^{n+1}+\frac{r}{\epsilon}\boldsymbol{H}^{n}, \\
\end{array}\right. \\
& \left\{\begin{array}{l}
r\frac{\phi_{21}^{n+1}}{\delta_{t}}+\left(rv_{r}^{n}\frac{\partial \phi^{n}}{\partial r}+rv_{z}^{n}\frac{\partial \phi^{n}}{\partial z}\right)=\mathcal{L}_{d} \left[\frac{\partial}{\partial r}\left(r\frac{\partial \mu_{21}^{n+1}}{\partial r}\right)+\frac{\partial}{\partial z}\left(r\frac{\partial \mu_{21}^{n+1}}{\partial z}\right)\right], \\
r\mu_{21}^{n+1}=-\epsilon \left[\frac{\partial}{\partial r}\left(r\frac{\partial \phi_{21}^{n+1}}{\partial r}\right)+\frac{\partial}{\partial z}\left(r\frac{\partial \phi_{21}^{n+1}}{\partial z}\right)\right]+\frac{rs}{\epsilon}\phi_{21}^{n+1}, \\
\end{array}\right. \\
& \left\{\begin{array}{l}
r\frac{\phi_{22}^{n+1}}{\delta_{t}}-\mathcal{L}_{d} \left[\frac{\partial}{\partial r}\left(r\frac{\partial \mu_{22}^{n+1}}{\partial r}\right)+\frac{\partial}{\partial z}\left(r\frac{\partial \mu_{22}^{n+1}}{\partial z}\right)\right]=0, \\
r\mu_{22}^{n+1}=-\epsilon \left[\frac{\partial}{\partial r}\left(r\frac{\partial \phi_{22}^{n+1}}{\partial r}\right)+\frac{\partial}{\partial z}\left(r\frac{\partial \phi_{22}^{n+1}}{\partial z}\right)\right]+\frac{rs}{\epsilon}\phi_{22}^{n+1}+\frac{r}{\epsilon}\boldsymbol{H}^{n}.
\end{array}\right.
&
\end{align}
$\boldsymbol{STEP 2:}$
\begin{equation}
    \begin{split}
&\mathcal{R}e\frac{r}{2}\frac{\rho^{n+1}+\rho^{n}}{\delta_{t}}\mathbf{u}_{1}^{n+1}-\nabla\cdot\left(r\eta^{n+1}D(\mathbf{u}_{1}^{n+1})\right)+2\frac{\eta^{n+1}}{r}(0,v^{n+1}_{1_{r}})=\mathcal{R}er\rho^{n}\frac{\tilde{\mathbf{u}}^{n+1}}{\delta_{t}} -\nabla\cdot\left(r\sqrt{\eta^{n+1}\eta^{n}}D(\mathbf{u}^{n})\right), \\ 
&\mathcal{R}e\frac{r}{2}\frac{\rho^{n+1}+\rho^{n}}{\delta_{t}}\mathbf{u}_{2}^{n+1}-\nabla\cdot\left(r\eta^{n+1}D(\mathbf{u}_{2}^{n+1})\right)+2\frac{\eta^{n+1}}{r}(0,v^{n+1}_{2_{r}})\\=&-\mathcal{R}e\rho^{n}\left(r\mathbf{u}^{n}\cdot\nabla\right) \mathbf{u}^{n}  -\frac{\mathcal{R}e}{2}\left(\nabla(r\rho^{n+1}\mathbf{u}^{n})\right)\mathbf{u}^{n}+\nabla\cdot\left(r\eta^{n+1}D(\mathbf{u}^{n})\right)-r\nabla(2p^{n}-p^{n-1})-\frac{1}{2}\nabla\cdot r\mathbf{J}(\mu^{n})\mathbf{u}^{n}-r\mathbf{J}(\mu^{n})\cdot\nabla\mathbf{u}^{n}.
    \end{split}
\end{equation}
$\boldsymbol{STEP 3:}$
\begin{align}
& p_{1}^{n+1}=p^{n}\\
& \frac{\partial}{\partial r}\left(r\frac{\partial p_{2}^{n+1}}{\partial r}\right)+\frac{\partial}{\partial z}\left(r\frac{\partial p_{2}^{n+1}}{\partial z}\right)=\frac{\chi\mathcal{R}e}{\delta_{t}}\left(\frac{\partial (rv_{r}^{n+1})}{\partial r}+\frac{\partial (rv_{z}^{n+1})}{\partial z}\right).
\end{align}

It can be observed that the overall computational cost of solving the numerical scheme \eqref{31}-\eqref{310} at each time step involves three elliptic systems with constant coefficients, two elliptic equations with positive variable coefficients, and one pressure Poisson equation with constant coefficients. Importantly, all these equations are fully decoupled, resulting in highly efficient calculations in practice.
\section{Numerical results}

In this section, we investigate the accuracy, energy stability, and effectiveness of the proposed scheme \eqref{31}-\eqref{310} numerically. We perform several numerical simulations to confirm the convergence rate and energy stability and compare our numerical results with physical experiments in \textcolor{red}{\cite{ref23}}. We also examine how the dynamics of droplet formation depend on various physical parameters of the system. To conduct the simulations, we use a rectangular region $[0,20] \times[0,3]\in \mathbb{R}^2$ as the computational domain and discretize space using the Finite Element Method (FEM) with a mesh size of 200$\times$30 unless otherwise specified. 
It should be noted that for the convenience of display, all the droplet formation diagrams are drawn vertically, and the coordinate axis in this section is rotated 90 degrees clockwise with respect to Figure \ref{fig:Fig 1}.
\subsection{Accuracy and stability test}

\begin{figure}[H]
\centering 
\begin{minipage}[b]{0.49\textwidth} 
\centering 
\includegraphics[width=55mm]{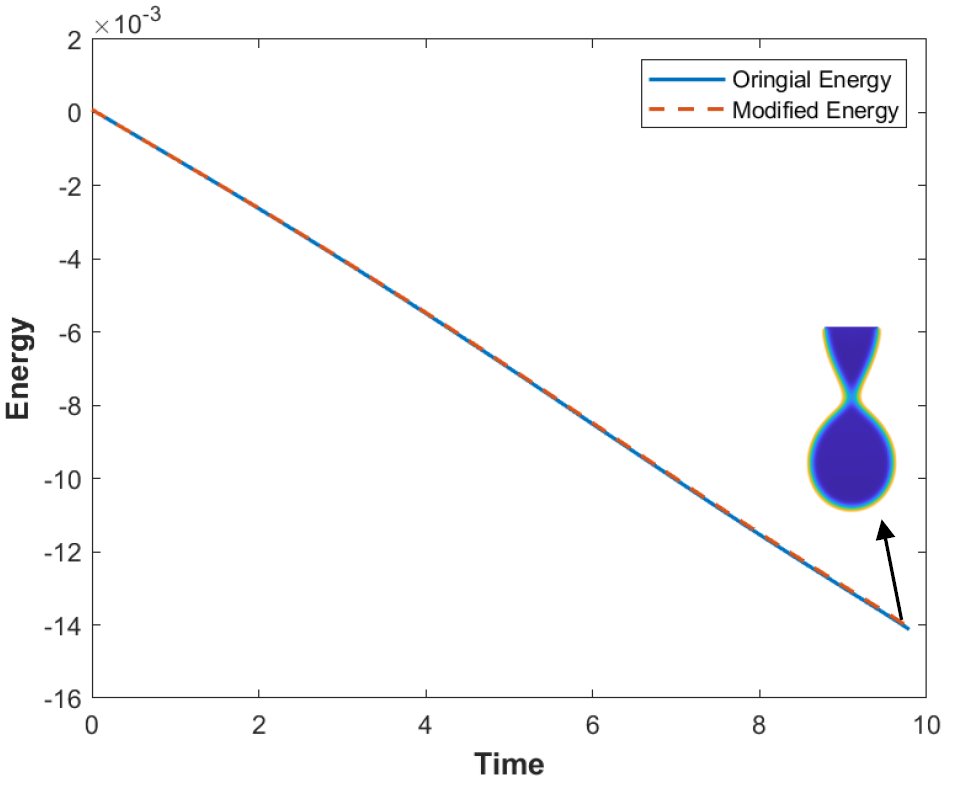} 
\caption{Energy evolutions.}
\label{fig:Fig.2}
\end{minipage}
\begin{minipage}[b]{0.49\textwidth} 
\centering 
\includegraphics[width=55mm]{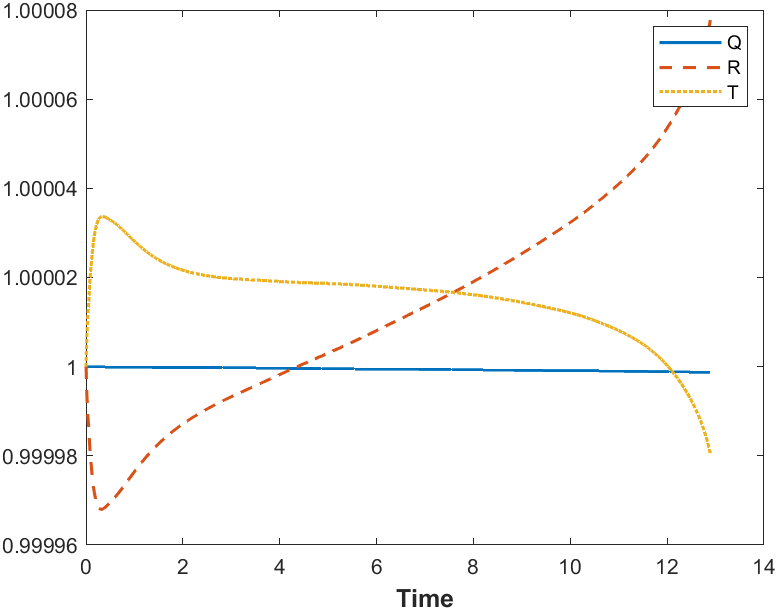} 
\caption{Q, R, T's evolutions.}
\label{fig:Fig.3}
\end{minipage}
\end{figure}

\begin{figure}[H]
\centering 
\begin{minipage}[b]{0.5\textwidth} 
\centering 
\includegraphics[width=60mm]{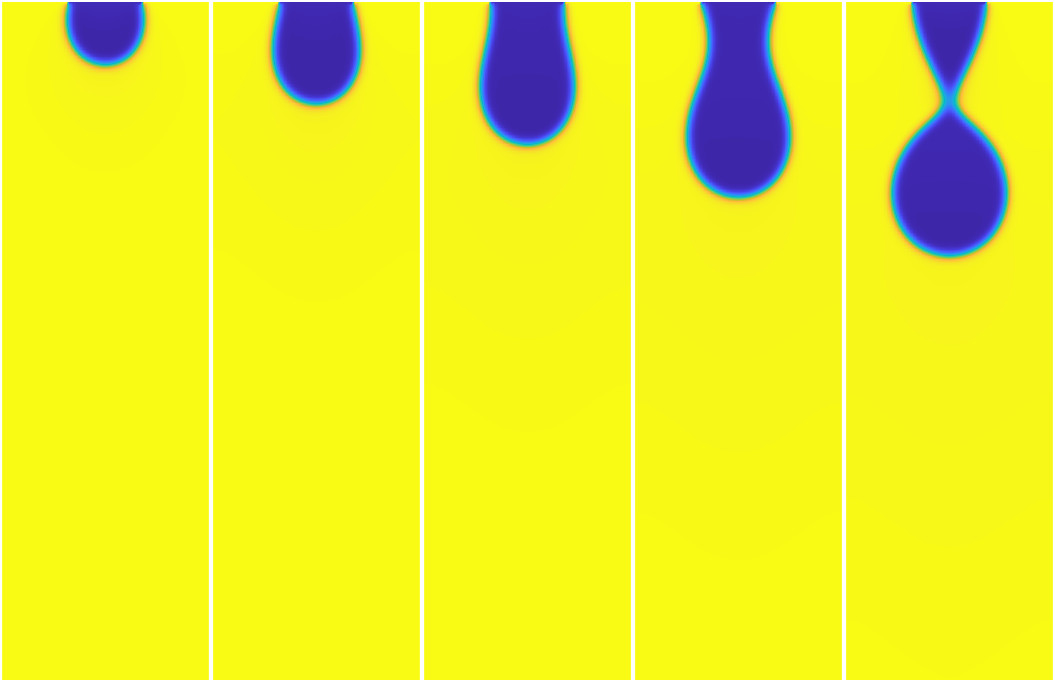} 
\caption{The process of droplet formation.}
\label{fig:Fig.4}
\end{minipage}
\begin{minipage}[b]{0.45\textwidth} 
\centering 
\includegraphics[width=55mm]{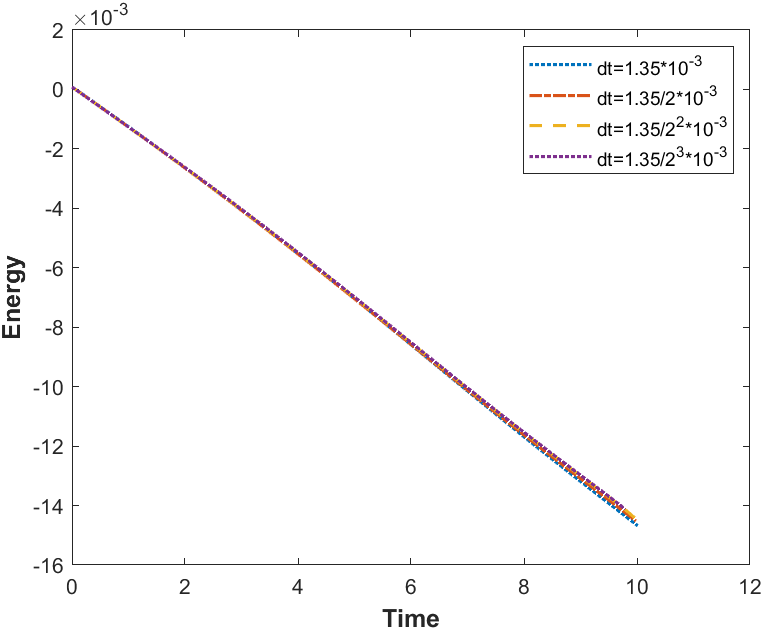} 
\caption{Stability tests.}
\label{fig:Fig.5}
\end{minipage}
\end{figure}

We perform an accuracy and stability test in this subsection. The model parameters are set as $a=3$, $S_{l}=20$, $n_z=200$, $n_r=30$, $\epsilon=0.1$, $\delta t=1.37\times 10^{-3}$, $Q_{r}=10$, $\alpha=10^{-3}$, $\mathcal{L}_{d}=0.05$, $\lambda_{\rho}=10$, $\lambda_{\eta}=1$,$Re=0.01$, $C a=0.04$, $\mathcal{B}=3/(2 \sqrt{2} \mathrm{Ca})$. Fig.\ref{fig:Fig.2} shows the comparison of original energy $\mathbf{E}^{n+1}_{O}$ in \eqref{223} and the modified energy $\mathbf{E}^{n+1}_{M}$ in \eqref{EM} after substracting a constant which we have mentioned in \eqref{dE} during the process of droplet formation. The states of the droplet at the endpoints of energy are also provided. We found that the two energy curves are very consistent, demonstrating the effectiveness of our algorithm. Fig.\ref{fig:Fig.3} shows the variations of the auxiliary variables $\boldsymbol{Q}$, $\boldsymbol{R}$ and $\boldsymbol{T}$ during the process of droplet formation. We found that their values are stable around 1, which further supports the effectiveness of our algorithm. As previously mentioned, the original system and the modified system are equivalent only when the values of these three auxiliary variables are around 1. Fig.\ref{fig:Fig.4} presents images of the droplet falling process at 1.67, 2.67, 5.34, 8.01, 12.89. Fig.\ref{fig:Fig.5} shows the evolution in energy for different time step sizes of $\delta_{t}=1.37\times 10^{-3}$, $1.37/2\times 10^{-3}$, $1.37/2^{2}\times 10^{-3}$, $1.37/2^{3}\times 10^{-3}$. The energy curves match well for different time step sizes, validating the convergence of our algorithm.

By varying the spatial resolution in this case, we can observe the convergence behavior of the numerical solution. Specifically, we consider the numerical solutions obtained with a fine grid, namely $h=1/160$, as an approximation of the exact solution. We consistently set the step size as $h=3/n_r=20/n_z$, where $n_r$ and $n_z$ represent the number of elements in the discretization of the $r$ and $z$ directions, respectively. The error with respect to the spatial step size $h$ is defined as  
\[\|e(\cdot)\|_{L^2}=\left[\sum_{i}\left|(\cdot)_{\text{exact},i}-(\cdot)_{h,i}\right|^2 \mathrm{~d} V_i\right]^{1 / 2},
\]
where $i$ iterates over all grid points, $(\cdot)_{\text {exact},i}$ are the values of the exact solution, $(\cdot)_{h,i}$ are the values of the numerical solution with a spatial step size h at the $i-$th grid point, and $\mathrm{~d} V_i$ represents the volume element associated with the grid point. Tabel~\ref{err_space} demonstrate that the chosen spatial discretization scheme exhibits a consistent and expected 2nd order, affirming its reliability in approximating the spatial aspects of the underlying problem.

\begin{table}
\begin{center}
\begin{tabular}{|c|c|c|c|c|c|c|c|c|}
\hline $h$ & $\|e(\phi)\|_{L^2}$ & Order & $\|e(v_z)\|_{L^2}$ & Order & $\|e(v_r)\|_{L^2}$ & Order & $\|e(p)\|_{L^2}$ & Order \\
\hline $1/10$ & $ 7.04  \mathrm{e}-01$ & & $2.73 \mathrm{e}-01$ & & $ 1.28 \mathrm{e}-01$ & & $1.06 \mathrm{e}+00$ & \\
\hline $1/20$ & $6.71 \mathrm{e}-02$ & 3.39 &  $1.57 \mathrm{e}-02$ & 4.11 & $6.90 \mathrm{e}-03$ & 4.21 & $1.74 \mathrm{e}-01$ & 2.61 \\
\hline $1/40$ & $1.46 \mathrm{e}-02$ & 2.20 & $3.90 \mathrm{e}-03$ & 2.01 & $1.70 \mathrm{e}-03$ & 2.02 & $4.18 \mathrm{e}-02$ & 2.05 \\
\hline $1/80$ & $3.60 \mathrm{e}-03$ & 2.02 & $  9.89\mathrm{e}-04$ & 1.98 & $4.36 \mathrm{e}-04$ & 1.96 & $1.08 \mathrm{e}-02$ & 1.95 \\
\hline 
\end{tabular}
\caption{\label{err_space}Grid refinement analysis with a fixed time step $\delta t = 10^{-5}$ at t=0.2. }
\end{center}
\end{table}

\begin{figure}[H]
\centering 
\begin{minipage}[b]{0.5\textwidth} 
\centering 
\includegraphics[width=60mm]{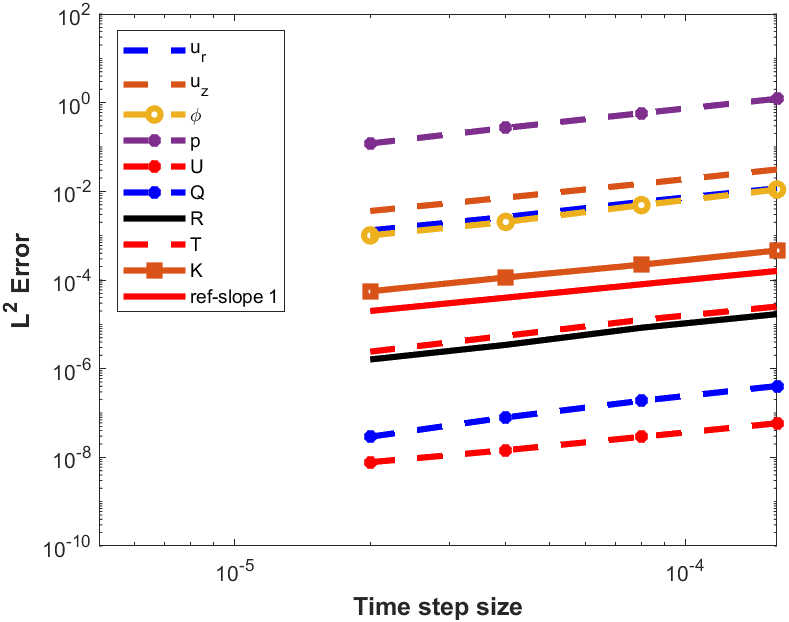} 
\caption{Accuracy tests.}
\label{fig:Fig.16}
\end{minipage}
\begin{minipage}[b]{0.45\textwidth} 
\centering 
\includegraphics[width=60mm]{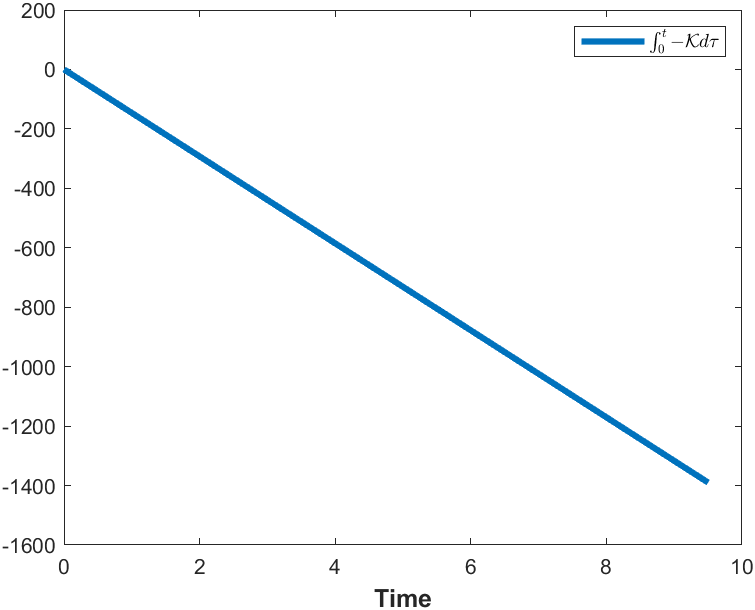} 
\caption{$\int_{0}^{t}-\mathcal{K}d\tau $'s evolution.}
\label{fig:Fig.17}
\end{minipage}
\end{figure}

Investigating the convergence behavior with respect to time, we conduct simulations with various time step sizes. Utilizing a sufficiently fine spatial mesh $800\times 120$, we compute the accuracy order of the temporal error relative to the time step size. For this analysis, we consider the numerical solutions obtained with a very small time step size, specifically $\delta_{t}=10^{-7}$, as the exact solution approximately. The $L^{2}$ errors for variables $\phi$, $\mathbf{u}$, $p$, $\boldsymbol{Q}$, $\boldsymbol{R}$, $\boldsymbol{T}$, $\boldsymbol{U}$, and $\boldsymbol{K}$ between the numerical solutions and the reference solutions at time $t = 0.2$ were calculated for different time step sizes. It should be noted that we know that the exact solutions of $\boldsymbol{Q}$, $\boldsymbol{R}$, and $\boldsymbol{T}$ are 1. The results, presented in Fig.\ref{fig:Fig.16}, demonstrate that all variables exhibit a very good first-order temporal convergence.

Fig.\ref{fig:Fig.17} shows the evolution of $\int_{0}^{t}-\mathcal{K}d\tau $ during droplet formation. As droplet formation is a finite time process, it is noteworthy that although the value of $\int_{0}^{t}-\mathcal{K}d\tau $ keeps decreasing, it remains bounded by a finite value throughout this process. So we only need to choose a sufficiently large $G$ to ensure that the definition of $\mathbf{K}$ is meaningful.

\subsection{Comparison with the experimental results by Zhang\cite{ref23}}

\begin{figure}[H]
\centering 
\begin{minipage}[b]{0.49\textwidth} 
\centering 
\includegraphics[width=50mm]{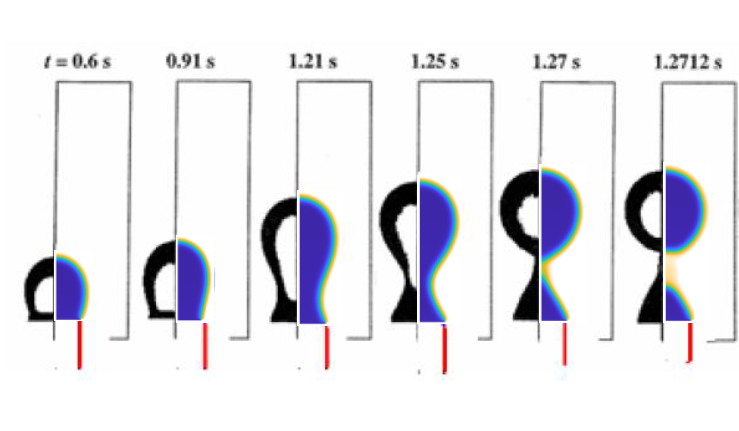} 
\caption{Comparison of the time sequences of bubble shapes.}
\label{fig:Fig.6}
\end{minipage}
\begin{minipage}[b]{0.49\textwidth} 
\centering 
\includegraphics[width=50mm]{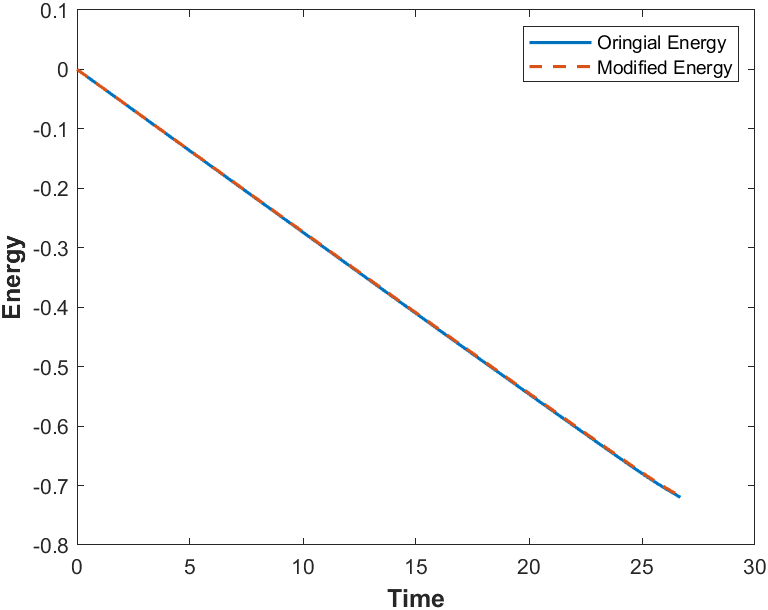}
\caption{Energy evolutions.}
\label{fig:Fig.7}
\end{minipage}
\end{figure}

Based on the CSF and VOF method, Zhang \cite{ref23} investigated bubble formation dynamics when a viscous liquid is injected through a vertical tube into another immiscible and viscous fluid. They focused on the dripping region where the dispersed fluid flowed through the capillary tube at a small flow rate. Good agreement was found between the numerical results and their experimental results. We now compare our numerical results with their experiment. The 2-ethyl-1-hexanol is the dispersed fluid and distilled water is the outer fluid. The viscosities of the inner and outer fluids are 0.089 g/(cm$\cdot$s) and 0.01 g/(cm$\cdot$s), the densities of the inner and outer fluids are 0.83 g/cm$^{3}$ and 1.0 g/cm$^{3}$, respectively. In this numerical case, We need to add the gravity term $\rho g$ to the momentum equation\eqref{23}. The other parameters are set as $a=4$, $S_{l}=20$, $n_z=200$, $n_r=40$, $\epsilon=0.1$, $Q_{r}=0$, $\mathcal{L}_{d}=0.156$, $\lambda_{\rho}=1.2048$, $\lambda_{\eta}=0.1123$, $Re=1.5461$, $C a=0.006986$, $\mathcal{B}=\frac{3}{2 \sqrt{2} \mathrm{Ca}}$, $\mathcal{B} o=1.5776$, $ \delta t=2.67\times 10^{-4}$, $\alpha=10^{-5}$, $g$=10.   Fig.\ref{fig:Fig.6} shows that the numerical results match well with the experimental results, thus validating the accuracy of our algorithm. As shown in Fig.\ref{fig:Fig.7}, the modified energy curve closely matches the original energy curve, thereby demonstrating the effectiveness of the algorithm.

\subsection{Effects of dimensionless parameters}

In this subsection, we explore the effects of various dimensionless parameters on the radius $R_d$ of the droplet, where $R_d$ is the maximum radius of the droplet before detaching from the nozzle. We assume the default values of $a=3$, $S_{l}=20$, $n_z=200$, $n_r=30$, $\epsilon=0.1$, $Q_{r}=10$. The time step size $ \delta t=2.67\times 10^{-3}$, and it is adjusted to a smaller value in some examples. 

\begin{figure}[h!]
\centering 
\begin{minipage}[b]{0.49\textwidth} 
\centering 
\includegraphics[width=50mm]{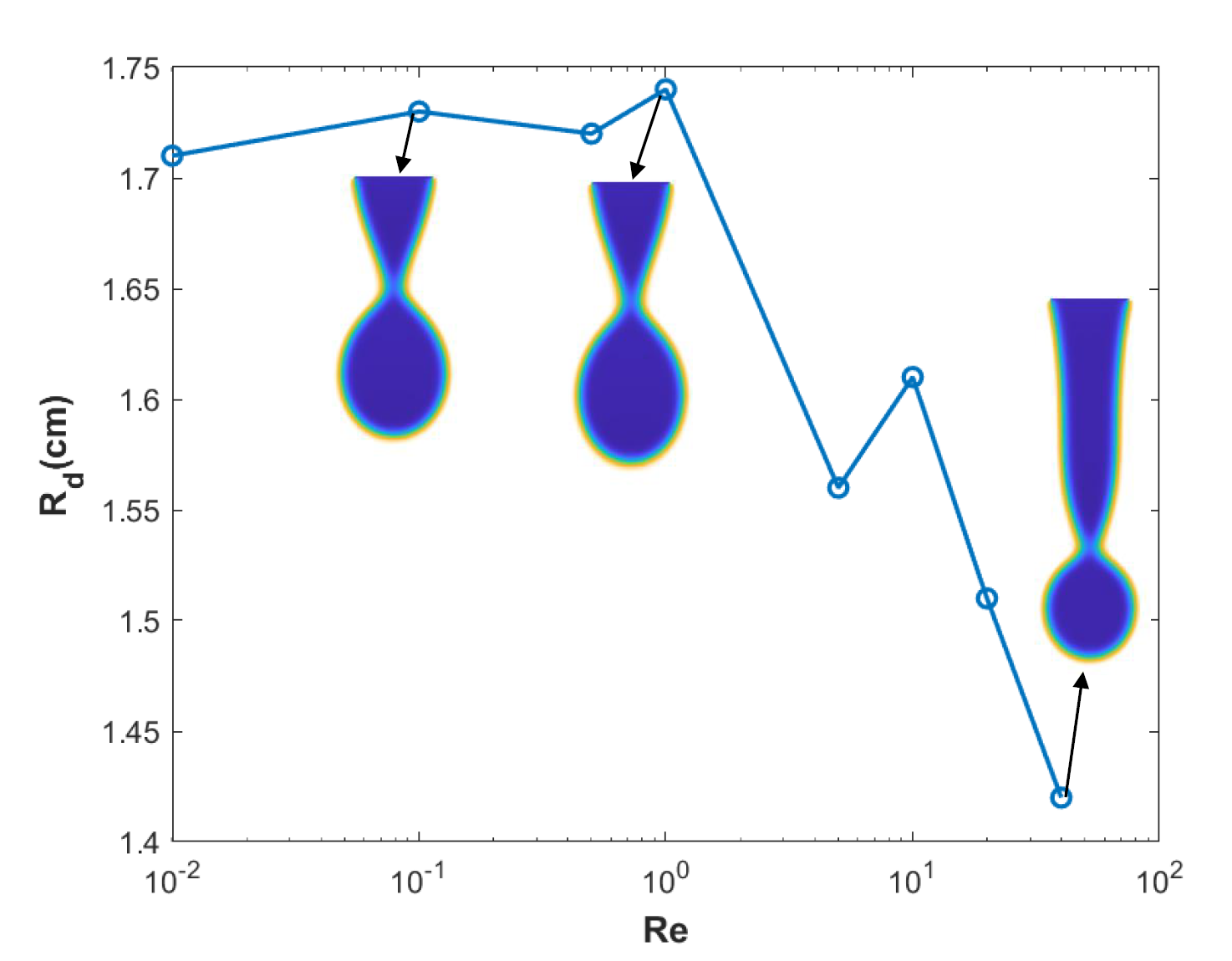} 
\caption{The effects of $\mathcal{R}$e on $R_{d}$.}
\label{fig:Fig.8}
\end{minipage}
\begin{minipage}[b]{0.49\textwidth} 
\centering 
\includegraphics[width=50mm]{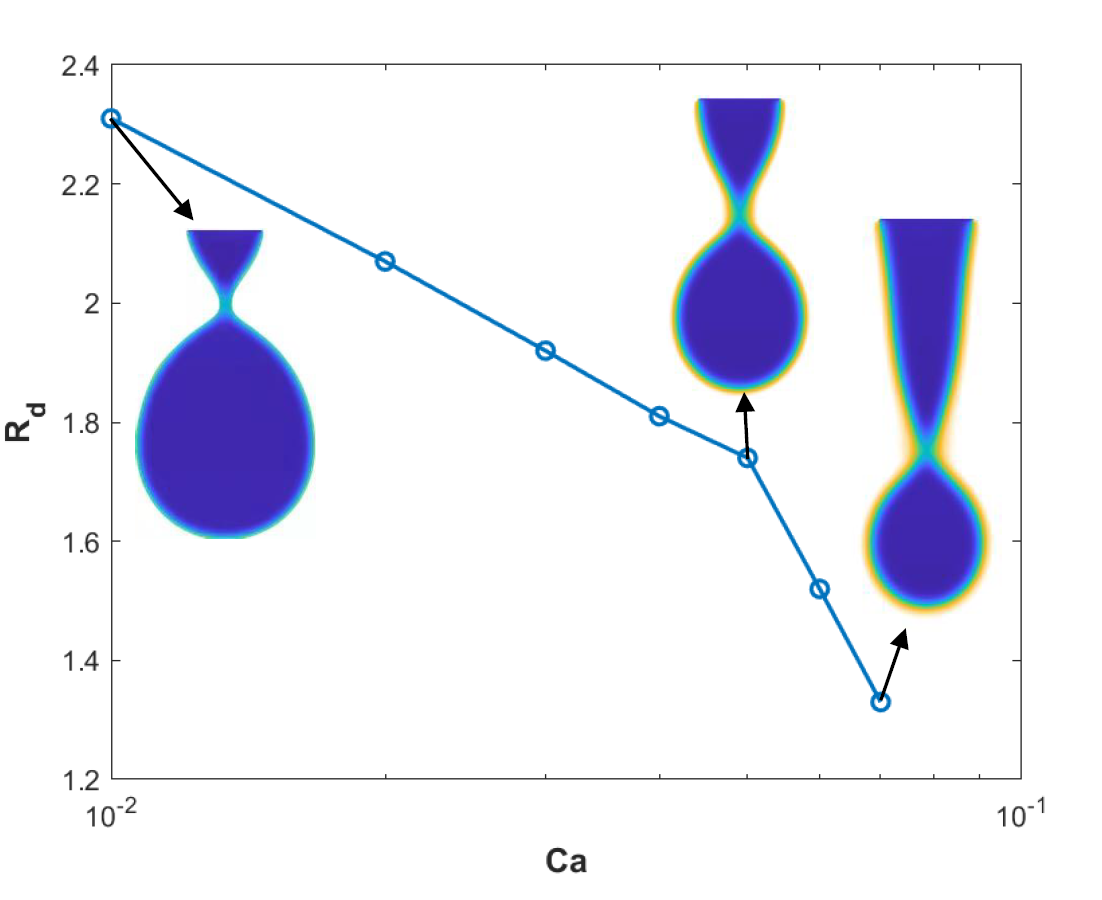}
\caption{The effects of $\mathcal{C}$a on $R_{d}$.}
\label{fig:Fig.9}
\end{minipage}
\begin{minipage}[b]{0.49\textwidth} 
\centering 
\includegraphics[width=50mm]{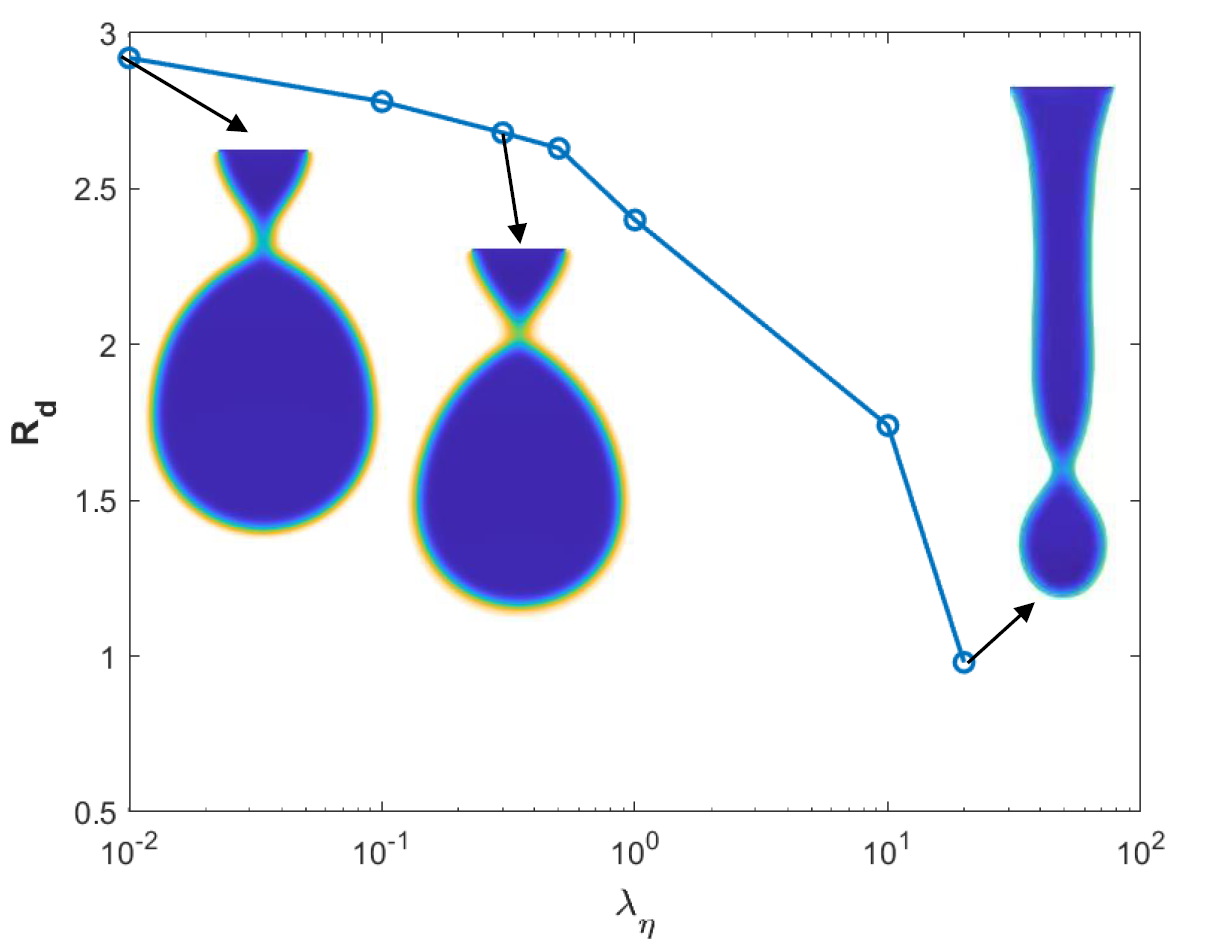} 
\caption{The effects of $\lambda_{\eta}$ on $R_{d}$.}
\label{fig:Fig.10}
\end{minipage}
\begin{minipage}[b]{0.49\textwidth} 
\centering 
\includegraphics[width=50mm]{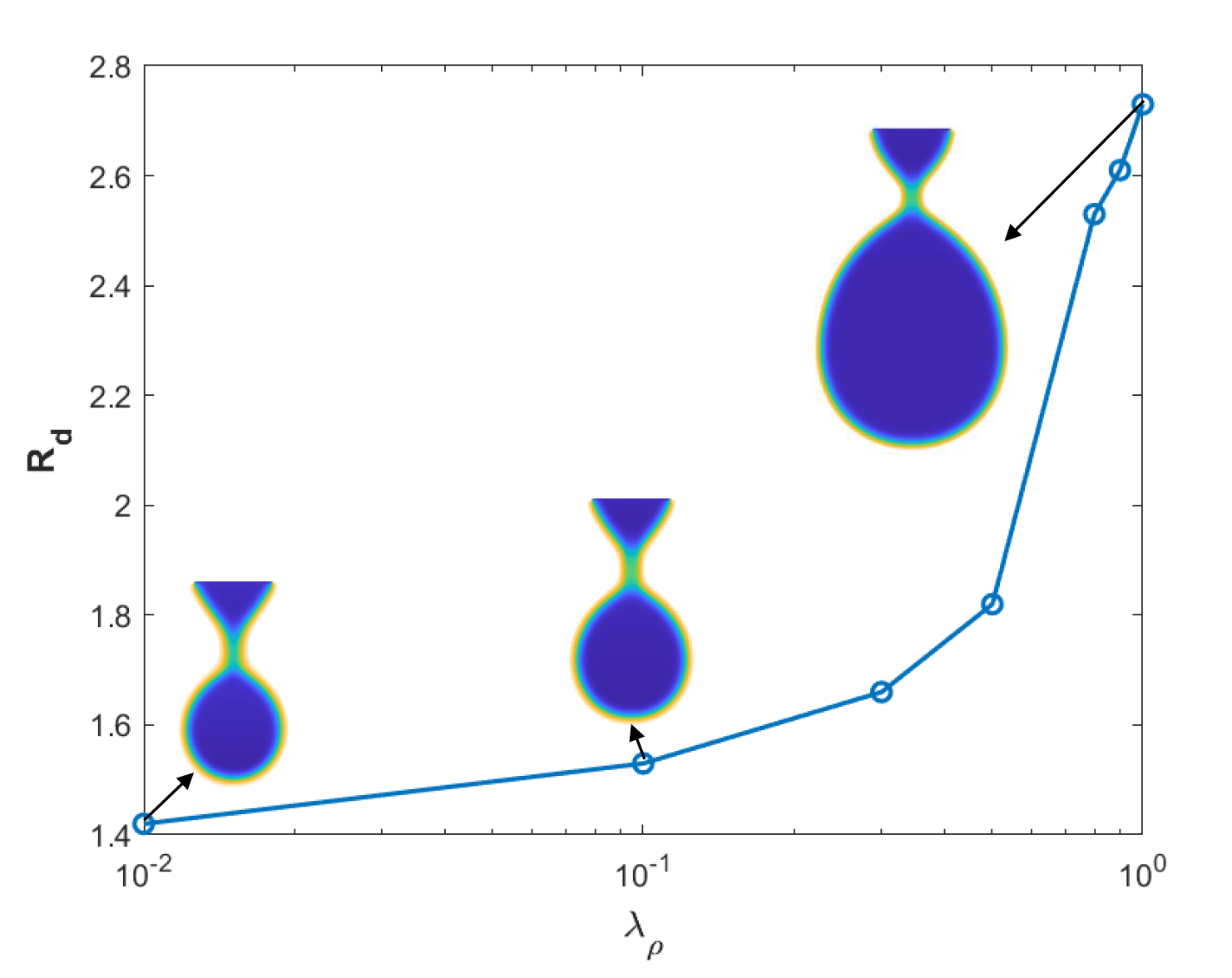}
\caption{The effects of $\lambda_{\rho}$ on $R_{d}$.}
\label{fig:Fig.11}
\end{minipage}
\end{figure}

\begin{figure}[H]
\centering 
\begin{minipage}[b]{0.49\textwidth} 
\centering 
\includegraphics[width=50mm]{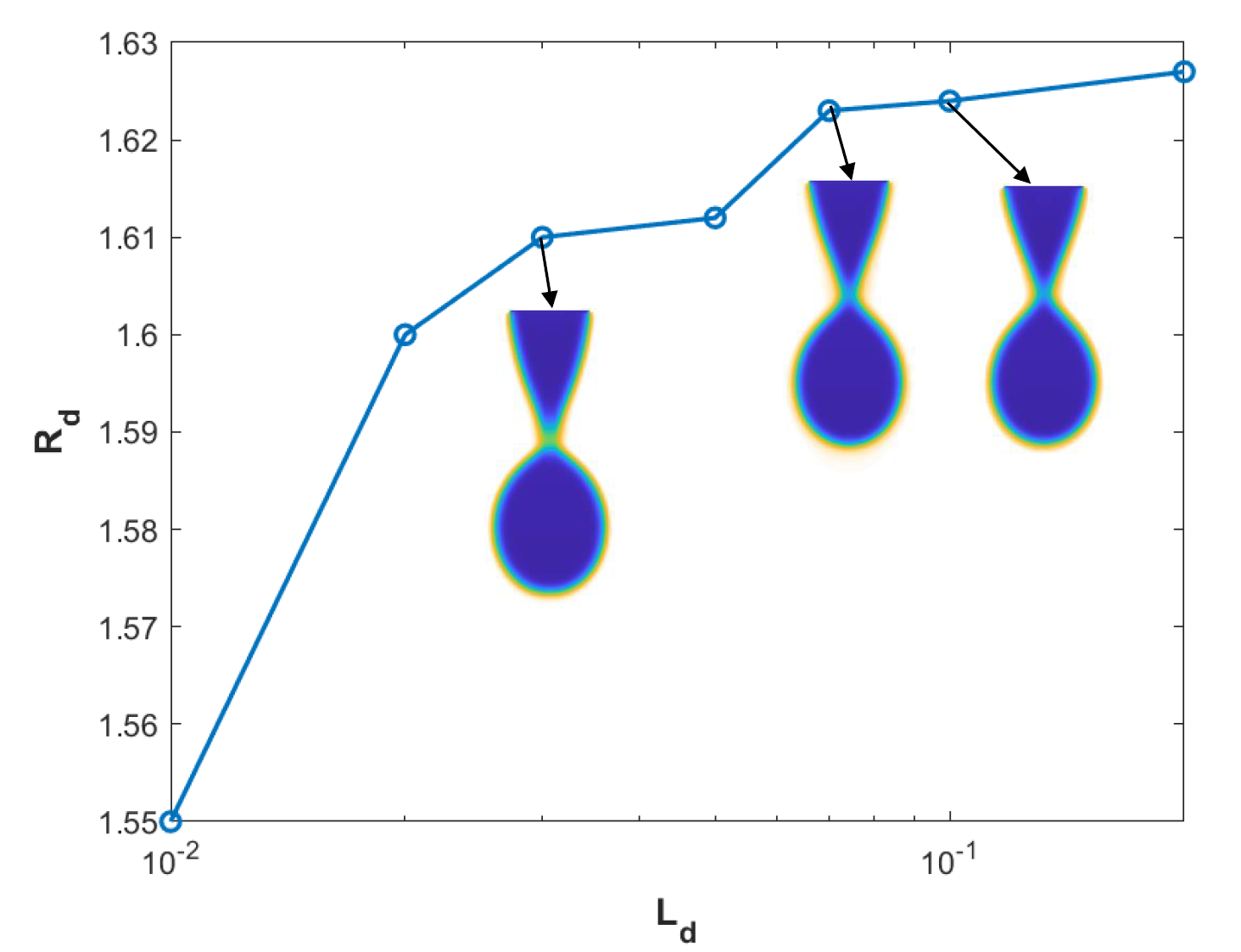} 
\caption{The effects of  $\mathcal{L}_{d}$ on $R_{d}$}
\label{fig:Fig.12}
\end{minipage}
\begin{minipage}[b]{0.49\textwidth} 
\centering 
\includegraphics[width=50mm]{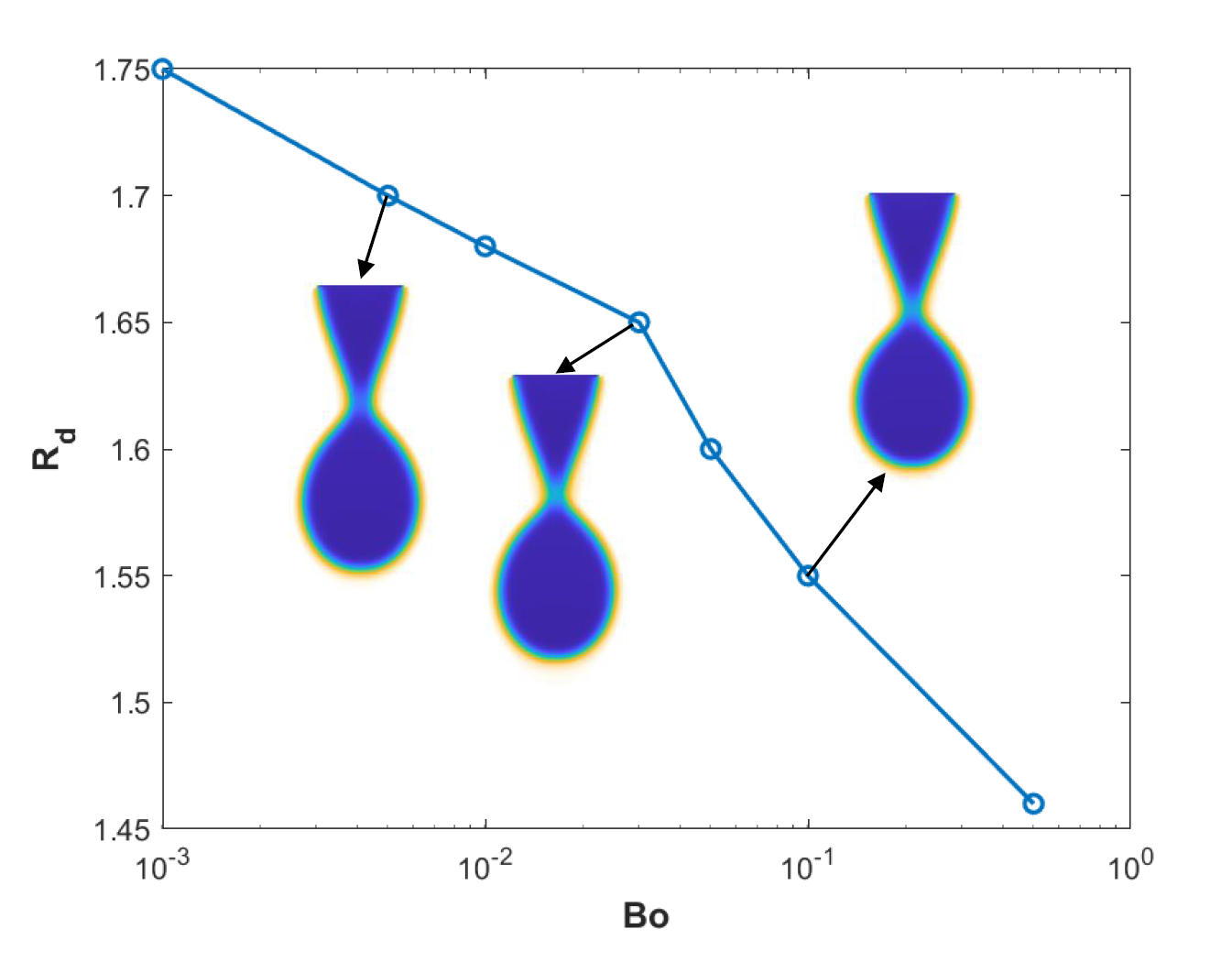}
\caption{The effects of $\mathcal{B}$o on $R_{d}$}
\label{fig:Fig.13}
\end{minipage}
\end{figure}

Fig.\ref{fig:Fig.8} shows the effects of  Reynolds number $\mathcal{R}$e on the maximum droplet radius $R_{d}$, where $\mathcal{R}$e is varied from left to right with values of 0.1, 1, and 40, and $\alpha$ is taken as $10^{-3}$, $10^{-3}$, $10^{-4}$, respectively. Other parameters are $\mathcal{L}_{d}=0.05$, $\lambda_{\rho}=0.1$, $\lambda_{\eta}=10$, $C a=0.01$ and $\mathcal{B}=\frac{3}{2 \sqrt{2} \mathrm{Ca}}$. By observing the relationship between $\mathcal{R}$e and $R_{d}$, we found that for $10^{-2}<\mathcal{R}e<10$, the influence of $\mathcal{R}$e on $R_{d}$ is small. However, as $\mathcal{R}$e increases for $\mathcal{R}e\ge 10$, $R_{d}$ decreases. Fig.\ref{fig:Fig.9} shows the effects of Capillary number $\mathcal{C}$a on the maximum droplet radius $R_{d}$, where $\mathcal{C}$a is varied from left to right with values of 0.01, 0.03, and 0.07. Other parameters are set as $\mathcal{L}_{d}=0.05$, $\lambda_{\rho}=0.1$, $\lambda_{\eta}=1$, $\mathcal{R}$e=0.01, $\mathcal{B}=\frac{3}{2 \sqrt{2} \mathrm{Ca}}$, $\alpha=10^{-4}$ and $ \delta t=2.67\times 10^{-4}$. By observing the relationship between $\mathcal{C}$a and $R_{d}$, we found that as $\mathcal{C}$a increases, $R_{d}$ decreases. Fig.\ref{fig:Fig.10} shows the effects of the viscosity ratio $\lambda_{\eta}$ on the maximum droplet radius $R_{d}$, where $\lambda_{\eta}$ is varied from left to right with values of 0.01, 0.5, and 20. $\alpha$ is set to $10^{-3}$, $10^{-4}$, $10^{-3}$ respectively. Other parameters are $\mathcal{L}_{d}=0.05$, $\lambda_{\rho}=0.1$, $\mathcal{R}$e=1, $\mathcal{C}$a=0.01, $\mathcal{B}=\frac{3}{2 \sqrt{2} \mathrm{Ca}}$. By observing the relationship between $\lambda_{\eta}$ and $R_{d}$, we found that as $\lambda_{\eta}$ increases, $R_{d}$ decreases. Fig.\ref{fig:Fig.11} shows the effects of density ratio $\lambda_{\rho}$ on $R_{d}$, where $\lambda_{\rho}$ is varied from left to right with values of 0.5, 0.8, and 1.0. $\alpha$ is set to $10^{-4}$, $10^{-5}$, $10^{-5}$ respectively. Other parameters are $\mathcal{L}_{d}=0.05$, $\lambda_{\eta}=0.1$, $\mathcal{R}$e=0.1, $\mathcal{C}$a=0.01,  $\mathcal{B}=\frac{3}{2 \sqrt{2} \mathrm{Ca}}$, and $ t=2.67\times 10^{-4}$. By observing the relationship between $\lambda_{\rho}$ and $R_{d}$, we found that as $\lambda_{\rho}$ increases, $R_{d}$ also increases. Fig.\ref{fig:Fig.12} shows the effects of the characteristic length $\mathcal{L}_{d}$ on the maximum droplet radius $R_{d}$, where $\mathcal{L}_{d}$ is varied from left to right with values of 0.001, 0.01, 0.1, and 1. Other parameters are $\lambda_{\rho}=10$, $\lambda_{\eta}=10$, $\mathcal{B}=\frac{3}{2 \sqrt{2} \mathrm{Ca}}$, $\mathcal{C}a=0.01$, $\mathcal{R}$e=1, and $\alpha=10^{-3}$. By observing the relationship between $\mathcal{L}_{d}$ and $R_{d}$, we found that as $\mathcal{L}_{d}$ increases, $R_{d}$ also increases. To account for gravity in the momentum equation, we need to add a gravity term $\frac{\mathcal{B}o}{\mathcal{C}a}\rho(g,0)$ while keeping the other equations unchanged, and the discretization format the same as before. Fig.\ref{fig:Fig.13} shows the effects of Bond number $\mathcal{B}$o on the maximum droplet radius $R_{d}$, where $\mathcal{B}$o is varied from left to right with values of 0.001, 0.01, 0.1, and 1. Other parameters are $\mathcal{L}_{d}=0.05$, $\lambda_{\rho}=0.1$, $\lambda_{\eta}=1$, $\mathcal{B}=\frac{3}{2 \sqrt{2} \mathrm{Ca}}$, $\mathcal{C}a=0.01$, $\mathcal{R}$e=1, $\alpha=10^{-3}$ and $g$ is the gravitational constant. By observing the relationship between $\mathcal{B}$o and $R_{d}$, we found that as $\mathcal{B}$o increases, $R_{d}$ decreases.

\section{Conclusions}

This paper proposes an efficient and unconditionally energy-stable numerical method for studying the dynamics of droplet formation. The method extends the numerical scheme proposed by Chen \cite{ref22} to address the droplet formation problem with nonhomogeneous Dirichlet boundary conditions. By combining the decoupled scalar auxiliary variable method with the modified penalty method and the operator Strang splitting method, we obtain an easy-to-implement numerical scheme. 
Our method is fully decoupled, linear, and unconditionally energy stable, offering significant computational efficiency. The proposed method only requires solving a few linear equations in each step, most of which have constant coefficients. We provide detailed information on the actual realization, solvability, and stability of the method. To demonstrate the accuracy and stability of the method, we present various numerical simulations.  Our simulation results demonstrate that the process of drop formation
can be reasonably predicted by the phase field model we used.

\section*{Acknowledgments}
X.-P. Wang acknowledges support from the National Natural Science Foundation of China
(NSFC) (No. 12271461), the key project of NSFC (No. 12131010), Shenzhen Science and Technology Innovation Program (Grant: C10120230046), the Hong Kong
Research Grants Council GRF (grants 16308421) and the University
Development Fund from The Chinese University of Hong Kong, Shenzhen (UDF01002028).
 C. Zhang is partially supported by the NSF of China (No.12101251), Guangzhou City Basic and Applied Basic Research Fund (No.2023A04J0010), and the National Key R$\&$D Program of China (2021YFA1002900).

\bibliographystyle{model1-num-names}
\bibliography{refs}


\end{document}